\documentclass{amsart}
\usepackage{amsxtra}
\usepackage{amssymb}
\newtheorem{Thm}{Theorem}[section]
\newtheorem{Prop}[Thm]{Proposition}
\newtheorem{Lem}[Thm]{Lemma}
\newtheorem{Cor}[Thm]{Corollary}
\newtheorem{Def}[Thm]{Definition}

\theoremstyle{remark}
\newtheorem{rem}[Thm]{Remark}
\numberwithin{equation}{section}

\begin{document}
\title[Quantized flag manifolds and irreducible $*$-representations]{Quantized 
flag manifolds and\\irreducible $*$-representations}
\author{Jasper V. Stokman}
\address{Department of Mathematics, University of Amsterdam, 
Plantage Muidergracht 24, 1018 TV Amsterdam, The Netherlands.}
\email{jasper@wins.uva.nl}
\author{Mathijs S. Dijkhuizen}
\address{Department of Mathematics, Faculty of Science,
Kobe University, Rokko, Kobe 657, Japan}
\email{msdz@math.s.kobe-u.ac.jp}
\thanks{The first author received financial support by NWO/Nissan}
\date{February, 1998}
\begin{abstract}
We study irreducible $*$-representations of a certain
quantization of the algebra of polynomial functions on
a generalized flag manifold regarded as a real manifold. All irreducible
$*$-representations are classified for a subclass of 
flag manifolds containing in particular the irreducible compact 
Hermitian symmetric spaces. For this subclass it is shown that the 
irreducible $*$-representations are parametrized by the symplectic
leaves of the underlying Poisson bracket.
We also discuss the relation between the quantized 
flag manifolds studied in this paper and the quantum flag manifolds
studied by Soibel'man, Lakshimibai \& Reshetikhin, Jur\v co \& \v S\v tov\'\i\v cek
and Korogodsky. 
\end{abstract}
\maketitle
\section{Introduction}
\label{section:intro}
The irreducible $*$-representations
of the ``standard'' quantization ${\mathbb{C}}_q\lbrack U \rbrack$ of the algebra  
of functions on a compact connected simple Lie group $U$ were 
classified by Soibel'man \cite{S}. He showed 
that there is a 1--1 correspondence between the
equivalence classes of irreducible $*$-representations of ${\mathbb{C}}_q\lbrack U
\rbrack$ and the symplectic leaves of the underlying Poisson bracket on $U$
(cf.\ \cite{S0}, \cite{S}).

This Poisson bracket is sometimes called Bruhat-Poisson, because its
symplectic foliation is a refinement of the Bruhat decomposition of $U$
(cf.\ Soibel'man \cite{S0}, \cite{S}). 
The symplectic leaves are naturally parametrized 
by $W\times T$, where $T\subset U$
is a maximal torus and $W$ is the Weyl group associated with $(U,T)$.

The 1--1 correspondence between 
equivalence classes of irreducible $*$-rep\-res\-enta\-tions 
of ${\mathbb{C}}_q\lbrack U
\rbrack$ and symplectic leaves of $U$ can be formally explained by the
observation that in the semi-classical limit the 
kernel of an irreducible $*$-representation should tend to a maximal Poisson ideal. 
The quotient of the Poisson algebra of polynomial functions on $U$ by this ideal 
is isomorphic to the Poisson algebra of functions on the symplectic leaf. 

In recent years many people have studied quantum homogeneous spaces (see
for example \cite{VS}, \cite{NYM}, \cite{Podles}, \cite{K0}, \cite{Noumi}, 
\cite{PV}, \cite{DS}). The results referred to above raise the obvious question 
whether the irreducible $*$-representations of
quantized function algebras on $U$-homogeneous spaces can be classified and related
to the symplectic foliation of the underlying Poisson  bracket.
This question was already raised in a paper by Lu \& Weinstein 
\cite[Question 4.8]{LW}, where they studied certain Poisson brackets on
$U$-homogeneous spaces that arise as a quotient of the Bruhat-Poisson bracket on $U$.

To our knowledge, affirmative answers to the above mentioned question 
have been given so far for only three different types of $U$-homogeneous spaces, 
namely Podle\'s's family of quantum 2-spheres \cite{Podles} (the relation with
the symplectic foliation of certain covariant Poisson brackets on the
2-sphere seems to have been observed for the first time by Lu \& Weinstein
\cite{LW2}), odd-dimensional complex quantum spheres 
$SU(n+1)/SU(n)$ (cf.\ Vaksman \& Soibel'man
\cite{VS}), and Stiefel manifolds $U(n)/U(n-l)$ 
(cf.\ Podkolzin \& Vainerman \cite{PV}).

In this paper we study the irreducible $*$-representations
of a certain quantized $\ast$-algebra of functions on a generalized flag manifold.
To be more specific, let $G$ denote the complexification of $U$, and let 
$P\subset G$ be a parabolic subgroup containing the standard Borel subgroup
$B_+$ with respect to a fixed choice of Cartan subalgebra and system of positive
roots (compatible with the choice of Bruhat-Poisson bracket on $U$, see
\cite{LW}). The generalized flag manifold $U/K$ naturally
becomes a Poisson $U$-homogeneous space (cf. Lu \& Weinstein \cite{LW}). 
The quotient Poisson bracket on 
$U/K$ is also called Bruhat-Poisson in \cite{LW}, since its symplectic leaves
coincide with the Schubert cells of the flag manifold $G/P \simeq U/K$.

It is straightforward to realize a quantum analogue ${\mathbb{C}}_q\lbrack K\rbrack$
of the algebra of polynomial functions on $K$ as a quantum subgroup of
${\mathbb{C}}_q\lbrack U\rbrack$. The corresponding $\ast$-subalgebra 
${\mathbb{C}}_q\lbrack U/K\rbrack$ of ${\mathbb{C}}_q\lbrack K\rbrack$-invariant
functions in ${\mathbb{C}}_q\lbrack U\rbrack$ may be regarded as a quantization
of the Poisson algebra of functions on $U/K$ endowed with the Bruhat-Poisson bracket.
The main result in this paper is a classification of all
the irreducible $*$-representations of ${\mathbb{C}}_q\lbrack U/K\rbrack$
for an important subclass of flag manifolds containing
in particular the irreducible Hermitian symmetric spaces of compact type. 
For this subclass we show that the equivalence classes of irreducible
$*$-representations are parametrized by the Schubert cells of $U/K$.
Let us emphasize that we regard here the flag manifold $U/K$ as
a real manifold. This means that the algebra of functions on $U/K$
has a natural $*$-structure, which survives quantization and allows us
to study $*$-representations in a way analogous to Soibel'man's
approach \cite{S}.

For an arbitrary generalized flag manifold $U/K$ we 
describe in detail how irreducible $*$-representations 
of ${\mathbb{C}}_q\lbrack U\rbrack$  decompose under restriction 
to ${\mathbb{C}}_q\lbrack U/K \rbrack$. 
This decomposition corresponds precisely to the way symplectic leaves
in $U$ project to Schubert cells in the flag manifold $U/K$. It leads immediately
to a classification of the irreducible $*$-representations
of the $C^*$-algebra $C_q(U/K)$, where $C_q(U/K)$ is obtained by taking the closure of
${\mathbb{C}}_q\lbrack U/K \rbrack$ with respect to the universal $C^\ast$-norm
on ${\mathbb{C}}_q[U]$. The equivalence classes of irreducible $*$-representations 
of $C_q(U/K)$ are naturally parametrized by the symplectic leaves
of $U/K$ endowed with the Bruhat-Poisson bracket.

For the classification of the irreducible $*$-representations
of the quantized function algebra ${\mathbb{C}}_q\lbrack U/K \rbrack$ itself
it is important to have a kind of
Poincar{\'e}-Birkhoff-Witt (PBW) factorization 
of ${\mathbb{C}}_q\lbrack U/K \rbrack$ (which in turn
is closely related to the irreducible decomposition of tensor products
of certain finite-dimensional irreducible $U$-modules).
Such a factorization is needed in order to develop 
a kind of highest weight representation theory for 
${\mathbb{C}}_q\lbrack U/K \rbrack$. In Soibel'man's
paper \cite{S}, a crucial role is played
by a similar factorization of ${\mathbb{C}}_q\lbrack U \rbrack$.
{}From Soibel'man's results one easily derives a factorization of the
algebra ${\mathbb{C}}_q\lbrack U/T \rbrack$ 
(corresponding to $P$ minimal parabolic in $G$).

In this paper we derive a PBW type factorization for a
different subclass of flag manifolds using the so-called 
Parthasarathy-Ranga Rao-Varadarajan (PRV) conjecture.
This conjecture was formulated as a follow-up to certain
results in the paper \cite{PRV} and was
independently proved by Kumar \cite{Ku} and Mathieu \cite{Ma} (see also 
Littelmann \cite{Li}).
The subclass of flag manifolds $U/K$ we consider here
can be characterized by the two conditions that
$(U,K)$ is a Gel'fand pair and that the Dynkin diagram of $K$ can be
obtained from the Dynkin diagram of $U$ by deleting one node (cf.\ 
Koornwinder \cite{Koo}).
Note that the corresponding $P\subset G$ is always maximal parabolic.
These two conditions are satisfied for the
irreducible compact Hermitian symmetric pairs $(U,K)$.

Roughly speaking, the PBW factorization in the above mentioned cases
states that the quantized function
algebra ${\mathbb{C}}_q[U/K]$ coincides with the quantized
algebra of zero-weighted complex valued polynomials on $U/K$. The
quantized algebra of zero-weighted complex valued polynomials can be naturally
defined for arbitrary generalized flag manifold $U/K$. It is always a  
$*$-subalgebra of ${\mathbb{C}}_q\lbrack U/K \rbrack$ and invariant
under the ${\mathbb{C}}_q\lbrack U\rbrack$-coaction (we shall call it the 
factorized $*$-algebra associated with $U/K$). 
The factorized $*$-algebra is closely related to the quantized algebra 
of holomorphic polynomials on generalized flag manifolds studied by
Soibel'man \cite{S2}, Lakshmibai \& Reshetikhin \cite{LR1}, \cite{LR2}, 
and Jur\v co \& \v S\v tov\'\i\v cek \cite{Jur} (for the classical groups)
as well as to the function spaces considered recently by Korogodsky
\cite{Kor}.

In this paper we classify the irreducible 
$*$-representations of the factorized $*$-algebra associated with an 
arbitrary flag manifold $U/K$ and we show 
that the equivalence classes of irreducible $*$-representations 
are naturally parametrized by the symplectic leaves
of $U/K$ endowed with the Bruhat-Poisson bracket. 
In particular, we obtain a complete classification
of the irreducible $*$-representations of 
${\mathbb{C}}_q\lbrack U/K \rbrack$ whenever a PBW type factorization holds 
for ${\mathbb{C}}_q\lbrack U/K \rbrack$  (i.e., 
${\mathbb{C}}_q\lbrack U/K \rbrack$ is equal 
to its factorized $*$-algebra).

The paper is organized as follows. In section 2 we review the results by
Lu \& Weinstein \cite{LW} and Soibel'man \cite{S} concerning
the Bruhat-Poisson bracket on $U$ and the
quotient Poisson bracket on a flag manifold. 
In section 3 we recall some well-known results
on the ``standard'' quantization of the universal 
enveloping algebra of a simple complex Lie algebra and its
finite-dimensional representations. We also recall the construction of 
the corresponding quantized function algebra ${\mathbb{C}}_q\lbrack U \rbrack$ and
give some commutation relations between certain matrix coefficients of
irreducible corepresentations of ${\mathbb{C}}_q\lbrack U \rbrack$.
They will play a crucial role in the classification of the irreducible
$*$-representations of the factorized $*$-algebra. 

In section 4 we define the quantized algebra 
${\mathbb{C}}_q\lbrack U/K\rbrack$ of functions 
on a flag manifold $U/K$ and its associated factorized $*$-subalgebra. 
We prove that the factorized $*$-algebra is equal to
${\mathbb{C}}_q\lbrack U/K \rbrack$ for the subclass of flag manifolds referred
to above.

In section 5 we study the restriction of an arbitrary 
irreducible $*$-representation of ${\mathbb{C}}_q\lbrack U \rbrack$ to 
${\mathbb{C}}_q\lbrack U/K \rbrack$. We use here
Soibel'man's explicit realization of the irreducible
$*$-representations of ${\mathbb{C}}_q\lbrack U \rbrack$ as tensor products of 
irreducible $*$-representations of ${\mathbb{C}}_q\lbrack SU(2) \rbrack$ 
(cf.\ \cite{S}, see also \cite{Ko}, \cite{VS} for $SU(n)$).
As a corollary we obtain a complete classification of the irreducible
$*$-representations of the $C^*$-algebra $C_q(U/K)$.

Section 6 is devoted to the classification of the
irreducible $*$-representations of the factorized $*$-algebra associated with an
arbitrary flag manifold. The techniques in section 6 are 
similar to those used by Soibel'man \cite{S}
for the classification of the irreducible $*$-representations of 
${\mathbb{C}}_q\lbrack U \rbrack$, and to those used by Joseph \cite{J}
to handle the more general problem of determining the 
primitive ideals of ${\mathbb{C}}_q\lbrack U\rbrack$.
\section{Bruhat-Poisson brackets on flag manifolds}
\label{section:Bruhat-Poisson}
In this section we review some results by Soibel'man \cite{S} and Lu \& Weinstein 
\cite{LW} concerning the Bruhat-Poisson bracket on a compact connected
simple Lie group $U$ and its flag manifolds.
For unexplained terminology in this section we refer 
the reader to \cite{CP} and \cite{LW}.
 
Let ${\mathfrak{g}}$ be a complex simple Lie algebra 
with a fixed Cartan subalgebra
${\mathfrak{h}}\subset {\mathfrak{g}}$. Let $G$ be the connected simply
connected Lie group with Lie algebra ${\mathfrak{g}}$ (regarded here as 
a real analytic Lie group).

Let $R\subset {\mathfrak{h}}^*$ be the root system associated with
$({\mathfrak{g}},{\mathfrak{h}})$ and write ${\mathfrak{g}}_{\alpha}$
for the root space associated with $\alpha\in R$. Let $\Delta=\lbrace
\alpha_1,\ldots,\alpha_r\rbrace$ be a basis of simple roots 
for $R$,  and let $R^+$ (resp. $R^-$) be the set of 
positive (resp. negative) roots relative to $\Delta$. We identify
${\mathfrak{h}}$ with its dual by the Killing form $\kappa$. 
The non-degenerate symmetric bilinear form on ${\mathfrak{h}}^*$ induced
by $\kappa$ is denoted by $(\cdot,\cdot )$. 
Let $W\subset \hbox{GL}({\mathfrak{h}}^*)$ be the
Weyl group of the root system $R$ and write $s_i=s_{\alpha_i}$ for the simple 
reflection associated with $\alpha_i\in\Delta$.

For $\alpha\in R$ write $d_{\alpha}:=(\alpha,\alpha)/2$.
Let $H_{\alpha}\in {\mathfrak{h}}$ be the element
associated with the coroot $\alpha\spcheck:=
d_{\alpha}^{-1}\alpha\in {\mathfrak{h}}^*$
under the identification $\mathfrak{h} \simeq \mathfrak{h}^*$.
Let us choose nonzero $X_{\alpha}\in {\mathfrak{g}}_{\alpha}$ ($\alpha\in R$)
such that for all $\alpha, \beta\in R$ one has
$\lbrack X_{\alpha},X_{-\alpha}\rbrack=H_{\alpha}$, 
$\kappa\bigl(X_{\alpha},X_{-\alpha}\bigr)=d_{\alpha}^{-1}$
and $\lbrack X_{\alpha},X_{\beta}\rbrack
=c_{\alpha,\beta}X_{\alpha+\beta}$ 
with $c_{\alpha,\beta}=-c_{-\alpha,-\beta}\in {\mathbb{R}}$
whenever $\alpha+\beta\in R$.
Let ${\mathfrak{h}}_0$ be the real form of ${\mathfrak{h}}$ defined as
the real span of the $H_{\alpha}$'s ($\alpha\in R$). Then
\begin{equation}\label{u}
{\mathfrak{u}}:=\sum_{\alpha\in R^+} \mathbb{R} (X_{\alpha}-X_{-\alpha})\oplus
\sum_{\alpha\in R^+} \mathbb{R} i(X_{\alpha}+X_{-\alpha})\oplus
i{\mathfrak{h}}_0
\end{equation}
is a compact real form of ${\mathfrak{g}}$.

Set ${\mathfrak{b}}:={\mathfrak{h}}_0\oplus {\mathfrak{n}}_+$
with ${\mathfrak{n}}_+:=\sum_{\alpha\in R^+}^\oplus{\mathfrak{g}}_{\alpha}$.
Then, by the Iwasawa decomposition for $\mathfrak{g}$, the triple
$({\mathfrak{g}},{\mathfrak{u}},{\mathfrak{b}})$ is a Manin triple
with respect to the imaginary part of the Killing form $\kappa$
(cf.\ \cite[\S4]{LW}).
Hence ${\mathfrak{u}}$, ${\mathfrak{b}}$ and ${\mathfrak{g}}$ 
naturally become Lie bialgebras. The dual Lie algebra 
${\mathfrak{u}}^*$ is isomorphic to
${\mathfrak{b}}$, and ${\mathfrak{g}}$ 
may be identified with the classical double of ${\mathfrak{u}}$.
The cocommutator $\delta: {\mathfrak{g}}\rightarrow {\mathfrak{g}}\wedge
{\mathfrak{g}}$ of the Lie bialgebra ${\mathfrak{g}}$ is coboundary, i.e.,
\[\delta(X)=\bigl(\hbox{ad}_X\otimes 1+1\otimes \hbox{ad}_X\bigr)r,\]
with the classical $r$-matrix $r\in {\mathfrak{g}}\wedge {\mathfrak{g}}$
given by the following well-known skew solution of the Modified 
Classical Yang-Baxter Equation,
\begin{equation}\label{r}
r=i\sum_{\alpha\in R^+}d_{\alpha}\bigl(X_{-\alpha}\otimes X_{\alpha}-X_{\alpha}\otimes
X_{-\alpha}\bigr)\in {\mathfrak{u}}\wedge {\mathfrak{u}}.
\end{equation}
The cocommutator on ${\mathfrak{u}}$ coincides with the 
restriction of ${\delta}$ to ${\mathfrak{u}}$. 

The corresponding Sklyanin bracket on the connected subgroup $U\subset G$ with
Lie algebra ${\mathfrak{u}}$ has 
\begin{equation}
\Omega_g = l_g^{\otimes 2} r - r_g^{\otimes 2} r
\end{equation}
as its associated Poisson tensor.
Here $l_g$ resp.\ $r_g$ denote infinitesimal left resp.\ right translation.
This particular Sklyanin bracket is often called Bruhat-Poisson, 
since its symplectic foliation
is closely related to the Bruhat decomposition of $G$. Let us
explain this in more detail.

Let $B$ be the
connected subgroup of $G$ with Lie algebra ${\mathfrak{b}}$, let
$T\subset U$ be the maximal torus in $U$ with Lie algebra $i{\mathfrak{h}}_0$,
and set $B_+:=TB$. 
The analytic Weyl group $N_U(T)/T$, where $N_U(T)$ is the normalizer of $T$ in $U$,
is isomorphic to $W$. More explicitly, the isomorphism sends 
the simple reflection $s_i$ to
$\exp\bigl(\frac{\pi}{2}(X_{\alpha_i}-X_{-\alpha_i})\bigr)/T$.
The double $B_+$-cosets in $G$ are parametrized by the elements of $W$.
Hence one has the Bruhat decomposition
\[
G=\coprod_{w\in W}B_+wB_+.
\]
By \cite[Prop.\  1.2.3.6]{W} the Bruhat decomposition has 
the following refinement:
\begin{equation}\label{refined}
G=\coprod_{m\in N_U(T)}BmB.
\end{equation}
For $m\in N_U(T)$ we set $\Sigma_{m}:=U\cap BmB$.
Then $\Sigma_{m}\not=\emptyset$ for all
$m\in N_U(T)$, and we have the disjoint union
\begin{equation}\label{Udecom}
U=\coprod_{m\in N_U(T)}\Sigma_{m}.
\end{equation}
Now recall that multiplication $U\times B\rightarrow G$ is a global diffeomorphism
by the Iwasawa decomposition of $G$.
So for any $b\in B$ and $u\in U$ there exists a unique $u^{b}\in U$ 
such that $bu\in u^bB$. As is easily verified, the map
\begin{equation}\label{rightdressing}
U\times B\rightarrow U,\quad (u,b)\mapsto u^{b^{-1}}
\end{equation}
is a right action of $B$ on $U$, and the corresponding decomposition of $U$ into
$B$-orbits coincides with the decomposition \eqref{Udecom}.
On the other hand, if we regard $B$ as the Poisson-Lie group dual to
$U$, the action \eqref{rightdressing} becomes the 
right dressing action of the dual group on $U$ (cf.\ \cite[Thm.\ 3.14]{LW}). 
Since the orbits in $U$ under the right dressing action
are exactly the symplectic leaves of the Poisson bracket on $U$ (cf.\ 
\cite[Thm.\ 13]{Se}, \cite[Thm.\ 3.15]{LW}), 
it follows that \eqref{Udecom} coincides with 
the decomposition of $U$ into symplectic leaves (cf.\ \cite[Theorem 2.2]{S}).

Next, we recall some results by Lu \& Weinstein \cite{LW} concerning 
certain quotient Poisson brackets on generalized flag manifolds.  
Let $S\subset \Delta$ be a set of simple roots, and
let $P_S$ be the corresponding standard parabolic subgroup of $G$.
The Lie
algebra ${\mathfrak{p}}_S$ of $P_S$ is given by
\begin{equation}\label{pS}
{\mathfrak{p}}_S:={\mathfrak{h}}\oplus\bigoplus_{\alpha\in\Gamma_S}
{\mathfrak{g}}_{\alpha}
\end{equation}
with $\Gamma_S:=R^+\cup\lbrace 
\alpha\in R \, | \, \alpha\in\hbox{span}(S)\rbrace$.
Let ${\mathfrak{l}}_S$ be the Levi factor of ${\mathfrak{p}}_S$,
\begin{equation}\label{lS}
{\mathfrak{l}}_S:={\mathfrak{h}}\oplus\bigoplus_{\alpha\in\Gamma_S\cap (-\Gamma_S)}
{\mathfrak{g}}_{\alpha},
\end{equation}
and set ${\mathfrak{k}}_S:={\mathfrak{p}}_S\cap {\mathfrak{u}}={\mathfrak{l}}_S\cap
{\mathfrak{u}}$.
Then ${\mathfrak{k}}_S$ is a compact real form of ${\mathfrak{l}}_S$.
Set $K_S:=U\cap P_S\subset U$, then $K_S\subset U$ is a Poisson-Lie subgroup of $U$ 
with Lie algebra ${\mathfrak{k}}_S$ (cf.\ \cite[Thm.\ 4.7]{LW}).
Hence there is a unique Poisson bracket on $U/K_S$ such that the
natural projection $\pi:U\rightarrow U/K_S$ is a Poisson map.
This bracket is also called Bruhat-Poisson. It is covariant in the
sense that the natural left action $U\times U/K_S\to U/K_S$ 
is a Poisson map.

Let $W_S$ be the subgroup of $W$ generated by the simple reflections in $S$.
One has $P_S = B_+W_S B_+$ (cf.\ \cite[Thm.\ 1.2.1.1]{W}). From this one
easily deduces that the double cosets $B_+xP_S$ ($x\in G$) are parametrized
by the elements of $W/W_S$. Hence one has
the Schubert cell decomposition of $U/K_S\simeq G/P_S$:
\begin{equation}\label{Schubertcell}
U/K_S=\coprod_{{\overline{w}}\in W/W_S}X_{\overline{w}},
\quad X_{\overline{w}}:=(U\cap B_+wP_S)/K_S\simeq B_+w/P_S,
\end{equation}
where ${\overline{w}}\in W/W_S$ is the right $W_S$-coset in $W$ which contains $w$.

Now, by \cite[Prop.\ 4.5]{LW}, the subgroup $K_S$ is invariant under
the action of $B$, which implies that the $B$-action descends to
$U/K_S$. The orbits in $U/K_S$ coincide exactly with the Schubert cells.
By \cite[Thm.\ 4.6]{LW} the symplectic leaves of the Poisson manifold
$U/K_S$ are exactly the orbits under the $B$-action. We conclude
(cf.\ \cite[Thm.\ 4.7]{LW}):
\begin{Thm}\label{Schubertleaf}
The decomposition into symplectic leaves of the flag manifold 
$U/K_S$ endowed with the Bruhat-Poisson bracket
coincides with its decomposition into
Schubert cells.
\end{Thm}
Consider now the set of minimal coset representatives
\begin{equation}\label{mincoset}
W^S:=\lbrace w\in W \, | \, l(ws_{\alpha})>l(w) \quad \forall \alpha\in S\rbrace.
\end{equation}
$W^S$ is a complete set of coset representatives for $W/W_S$. 
Any element $w\in W$ can be uniquely written as a product $w=w_1w_2$ with
$w_1\in W^S$, $w_2\in W_S$. The elements of $W^S$ are minimal in the sense that
\begin{equation}\label{minimal}
l(w_1w_2)=l(w_1)+l(w_2),\quad  (w_1\in W^S, w_2\in W_S),
\end{equation}
where $l(w):=\#(R^+\cap wR^-)$ is the length function on $W$.

Observe that $\pi$ maps the symplectic leaf 
$\Sigma_{m}\subset U$ onto the symplectic
leaf $X_{\overline{w(m)}}\subset U/K_S$, where $w(m):=m/T\in W$.
We write $\pi_{m}: \Sigma_{m}\rightarrow X_{\overline{w(m)}}$ 
for the surjective Poisson map obtained by restricting $\pi$ to the 
symplectic leaf $\Sigma_{m}$.
The minimality condition \eqref{mincoset} 
translates to the following property of the map $\pi_{m}$. 
\begin{Prop}\label{semiclassical}
Let $m\in N_U(T)$. Then $\pi_{m}:\Sigma_{m}\rightarrow
X_{\overline{w(m)}}$ is a symplectic automorphism if and only if $w(m)\in W^S$.
\end{Prop}
\begin{proof}
For $w\in W$ set 
\[
{\mathfrak{n}}_w:=\bigoplus_{\alpha\in R^+\cap wR^-}{\mathfrak{g}}_{\alpha},
\quad N_w:=\hbox{exp}({\mathfrak{n}}_w).
\]
Observe that the complex dimension of $N_w$ is equal to $l(w)$. 
Write $\text{pr}_U: G\simeq U\times B\rightarrow U$ for the canonical projection.
It is well known that for $m\in N_U(T)$ and for $w\in W^S$ with representative
$m_w\in N_U(T)$, the maps
\begin{equation}
\begin{split}
\phi_{m}: N_{w(m)}&\rightarrow \Sigma_{m},\quad
n\mapsto \text{pr}_U(nm),\nonumber\\
\psi_{w}:N_w&\rightarrow X_{\overline{w}},\quad n\mapsto
\pi\bigl(\text{pr}_{U}(nm_w)\bigr)\nonumber
\end{split}
\end{equation}
are surjective diffeomorphisms (see for example \cite[Proposition 1.1 \& 5.1]{BGG}).
The map $\psi_w$ is independent of the choice of representative $m_w$
for $w$. It follows now from \eqref{minimal} by a dimension count that $\pi_m$ 
can only be a diffeomorphism if $w(m)\in W^S$.
On the other hand, if $m\in N_U(T)$ such that $w(m)\in W^S$, then
$\pi_{m}=\psi_{w(m)}\circ{\phi}_m^{-1}$ and hence $\pi_m$ is
a diffeomorphism.
\end{proof}
Soibel'man \cite{S} gave a description of the symplectic leaves $\Sigma_{m}$
($m\in N_U(T)$) as a product of two-dimensional leaves
which turns out to have a nice generalization to the quantized setting
(cf.\ section 5).
For $i\in [1,r]$, let $\gamma_i: SU(2)\hookrightarrow U$ be the 
embedding corresponding to the $i$th node of the Dynkin diagram of $U$.
After a possible renormalization of the Bruhat-Poisson structure on $SU(2)$,
$\gamma_i$ becomes an embedding of Poisson-Lie groups.
Recall that the two-dimensional leaves of $SU(2)$ are given by
\begin{gather*}
S_t:=\left\{ 
\begin{pmatrix} \alpha & \beta \\ -{\overline{\beta}} & {\overline{\alpha}}
\end{pmatrix} \in \left.\hbox{SU}(2) \, \right|
 \, \hbox{arg}(\beta)=\hbox{arg}(t) \right\}
\quad (t\in {\mathbb{T}})
\end{gather*}
where ${\mathbb{T}}\subset {\mathbb{C}}$ is the unit circle in the complex plane.
The restriction of the embedding $\gamma_i$ to $S_1\subset SU(2)$ 
is a symplectic automorphism
from $S_1$ onto the symplectic leaf $\Sigma_{m_i}\subset U$, where 
$m_i=\hbox{exp}\bigl(\frac{\pi}{2}(X_{\alpha_i}-X_{-\alpha_i})\bigr)$. 
Recall that $m_i\in N_U(T)$ is a 
representative of the simple reflection $s_i\in W$.

For arbitrary $m\in N_U(T)$ let 
$w(m)=s_{i_1}s_{i_2}\cdots s_{i_l}$ be a reduced expression for $w(m):=m/T\in W$,
and let $t_m\in T$ be the unique element such that 
$m=m_{i_1}m_{i_2}\cdots m_{i_r}t_m$. Note that $t_m$ depends on the choice
of reduced expression for $w(m)$.  
The map
\[(g_1,\ldots,g_l)\mapsto \gamma_{i_1}(g_1)\gamma_{i_2}(g_2)\cdots
\gamma_{i_l}(g_l)t_m
\]
defines a symplectic automorphism from $S_1^{\times l}$ 
onto the symplectic leaf 
$\Sigma_{m}\subset U$ (cf.\ \cite[\S2]{S}, \cite{St}). 
Note that the image of the map is independent of the choice of 
reduced expression for $w(m)$, although the map itself is not.

Combined with Proposition \ref{semiclassical} we now obtain the following
description of the symplectic leaves of the generalized flag manifold $U/K_S$.
\begin{Prop}\label{finestructure}
Let $m\in N_U(T)$ and set $w:=m/T\in W$. Let $w_1\in W^S$, $w_2\in
W_S$ be such that $w=w_1w_2$ and 
choose reduced expressions $w_1=s_{i_1}\cdots s_{i_p}$ and
$w_2=s_{i_{p+1}}\cdots s_{i_l}$. 
Then the map
\[
(g_1,g_2,\ldots,g_l)\mapsto \gamma_{i_1}(g_1)\gamma_{i_2}(g_2)\cdots
\gamma_{i_l}(g_l)/K_S
\]
is a surjective Poisson map from $S_1^{\times l}$ onto the Schubert cell
$X_{\overline{w}}$. It factorizes through the projection 
$\text{pr}\colon S_1^{\times l} = S_1^{\times p} \times S_1^{\times (l-p)}
\rightarrow S_1^{\times p}$.
The quotient map from $S_1^{\times p}$ onto $X_{\overline{w}}$ is a
symplectic automorphism. In particular, we have
\[
X_{\overline{w}}=
\bigl(\Sigma_{m_{i_1}}\Sigma_{m_{i_2}}\cdots \Sigma_{m_{i_p}}\bigr)/K_S.
\qed \]
\renewcommand{\qed}{}\end{Prop}
See Lu \cite{Lu} for more details in the case of the full flag manifold
($K_S=T$).
\section{Preliminaries on the quantized function algebra 
${\mathbb{C}}_q\lbrack U \rbrack$}
In this section we introduce some
notations which we will need throughout the remainder of this paper.
First, we recall the definition of the quantized universal enveloping
algebra associated with the simple complex Lie algebra ${\mathfrak{g}}$.
We use the notations introduced in the previous section.

Set $d_i:=d_{\alpha_i}$ and $H_i:=H_{\alpha_i}$ for $i\in [1,r]$.
Let $A=(a_{ij})$ be the Cartan matrix, i.e.\ 
$a_{ij}:=d_i^{-1}(\alpha_i,\alpha_j)$.
Note that $H_i\in {\mathfrak{h}}$ is the unique element such
that $\alpha_j(H_i)=a_{ij}$ for all $j$. 
The weight lattice is given by
\begin{equation}
P=\lbrace \lambda\in {\mathfrak{h}}^* \, | \,
\lambda(H_i)=(\lambda,\alpha_i\spcheck)\in {\mathbb{Z}} \,\,\,\,\,\forall i\rbrace.
\end{equation}
The fundamental weights $\varpi_{\alpha_i}=\varpi_i$ ($i\in [1,r]$)
are characterized by $\varpi_i(H_j)=(\varpi_i,\alpha_j\spcheck)=\delta_{ij}$ 
for all $j$. The set of dominant weights $P_+$ resp.\ 
regular dominant weights $P_{++}$ is equal to ${\mathbb{K}}\hbox{-span}\lbrace
\varpi_{\alpha}\rbrace_{\alpha\in\Delta}$ with ${\mathbb{K}}={\mathbb{Z}}_+$
resp.\ ${\mathbb{N}}$.

We fix $q\in (0,1)$. The quantized universal enveloping algebra 
$U_q({\mathfrak{g}})$ associated with the simple Lie algebra
${\mathfrak{g}}$ is the unital associative algebra
over ${\mathbb{C}}$ with generators $K_i^{\pm 1}$, $X_i^{\pm}$ ($i=[1,r]$)
and relations
\begin{equation}
\begin{split}
&K_iK_j=K_jK_i,\quad K_iK_i^{-1}=K_i^{-1}K_i=1\\
&K_iX_j^{\pm}K_i^{-1}=q_i^{\pm\alpha_j(H_i)}X_j^{\pm}\\
&X_i^+X_j^--X_j^-X_i^+=\delta_{ij}\frac{K_i-K_i^{-1}}{q_i-q_i^{-1}}\\
&\sum_{s=0}^{1-a_{ij}}(-1)^s\binom{1-a_{ij}}{s}_{q_i}(X_i^{\pm})^{1-a_{ij}-s}
X_j^{\pm}(X_i^{\pm})^s=0 \quad (i\not=j)
\end{split}
\end{equation}
where $q_i:=q^{d_i}$, 
\[\lbrack a \rbrack_q:=\frac{q^a-q^{-a}}{q-q^{-1}} \quad (a\in {\mathbb{N}}),
\quad [0]_q:=1,\]
$\lbrack a \rbrack_q!:=\lbrack a \rbrack_q\lbrack a-1\rbrack_q\ldots \lbrack 1
\rbrack_q$,  and 
\[\binom{a}{n}_q:=\frac{\lbrack a \rbrack_q!}
{\lbrack a-n\rbrack_q!\lbrack n \rbrack_q!}.\]
A Hopf algebra structure on $U_q(\mathfrak{g})$ is uniquely determined by
the formulas
\begin{equation}\label{commrelationsU}
\begin{split}
&\Delta(X_i^+)=X_i^+\otimes 1 + K_i\otimes X_i^+, \quad
\Delta(X_i^-)=X_i^-\otimes K_i^{-1}+1\otimes X_i^-,\\
&\Delta(K_i^{\pm 1})=K_i^{\pm 1}\otimes K_i^{\pm 1},\\
&S(K_i^{\pm 1})=K_i^{\mp 1},\quad 
S(X_i^+)=-K_i^{-1}X_i^+,\quad  S(X_i^-)=-X_i^-K_i,\\
&\varepsilon(K_i^{\pm 1})=1,\quad \varepsilon(X_i^{\pm})=0.
\end{split}
\end{equation}
In fact, $U_q({\mathfrak{g}})$ may be regarded as a quantization of 
the co-Poisson-Hopf algebra structure (cf.\ \cite[Ch.\ 6]{CP})
on $U(\mathfrak{g})$ induced by
the Lie bialgebra $({\mathfrak{g}}, -i \delta)$,
$\delta$ being the cocommutator of ${\mathfrak{g}}$ associated with the
$r$-matrix \eqref{r}.
$U_q(\mathfrak{g})$ becomes a Hopf $*$-algebra with $*$-structure on the
generators given by
\begin{equation}\label{star}
(K_i^{\pm 1})^*=K_i^{\pm 1},\quad
(X_i^+)^*=q_i^{-1}X_i^-K_i,\quad (X_i^-)^*=q_iK_i^{-1}X_i^+.
\end{equation}
In the classical limit $q\to 1$,
the $*$-structure becomes an involutive, conjugate-linear 
anti-automorphism of ${\mathfrak{g}}$ with $-1$ eigenspace 
equal to the compact real form ${\mathfrak{u}}$ defined in \eqref{u}. 

Let $U^{\pm}=U_q({\mathfrak{n}}_{\pm})$ 
be the subalgebra of $U_q(\mathfrak{g})$ generated by
$X_i^{\pm}$ ($i=[1,r]$) and write $U^0:=U_q(\mathfrak{h})$ for the commutative
subalgebra generated by $K_i^{\pm 1}$ ($i=[1,r]$).
Let us write $Q$ (resp.\ $Q^+$) for the integral (resp.\ positive integral) 
span of the positive roots. We have the direct sum decomposition
\[U^{\pm}=\bigoplus_{\alpha\in Q^+}U^{\pm}_{\pm \alpha},\]
where 
$U^{\pm}_{\pm \alpha}:=\lbrace \phi\in U^{\pm} \, | \,
K_i\phi K_i^{-1}=q_i^{\pm\alpha(H_i)}\phi\rbrace$. 
The Poincar{\'e}-Birkhoff-Witt 
Theorem for $U_q(\mathfrak{g})$ states that multiplication defines an
isomorphism of vector spaces
\[U^-\otimes U^0\otimes U^+ \to U_q(\mathfrak{g}).\]
In particular, $U_q(\mathfrak{g})$ is spanned by elements of the form
$b_{-\eta}K^{\alpha}a_{\zeta}$ where $b_{-\eta}\in U_{-\eta}^-, a_{\zeta}\in
U_{\zeta}^+$ ($\eta,\zeta\in Q_+$) and $\alpha\in Q$.
Here we used the notation 
$K^{\alpha}=K_1^{k_1}\cdots K_r^{k_r}$ if $\alpha=\sum_ik_i\alpha_i$.

For a left $U_q(\mathfrak{g})$-module $V$, we say that
$0\not=v\in V$ has weight $\mu\in {\mathfrak{h}}^*$ if 
$K_i\cdot v=q_i^{\mu(H_i)}v=q^{(\mu,\alpha_i)}v$  
for all $i$.
We write $V_{\mu}$ for the corresponding weight space.
Recall that a $P$-weighted 
finite-dimensional irreducible representation of
$U_q(\mathfrak{g})$ is a highest weight module $V=V(\lambda)$ with highest
weight $\lambda\in P_+$. If $v_{\lambda}\in V(\lambda)$ is a 
highest weight vector, we
have $V(\lambda)=\sum_{\alpha\in Q^{+}}^{\oplus}U_{-\alpha}^-v_{\lambda}$ by 
the PBW Theorem, hence the set of weights $P(\lambda)$ of $V(\lambda)$ is 
a subset of the weight lattice $P$ satisfying $\mu\leq\lambda$ 
for all $\mu\in P(\lambda)$. Here $\leq$ is the dominance order on $P$
(i.e.\ $\mu\leq\nu$ if $\nu-\mu\in Q^+$ and 
$\mu<\nu$ if $\mu\leq\nu$ and $\mu\not=\nu$).
 
We define irreducible finite dimensional $P$-weighted right 
$U_q({\mathfrak{g}})$-modules with respect to the opposite Borel subgroup. So 
the irreducible finite dimensional right 
$U_q({\mathfrak{g}})$-module $V(\lambda)$ with highest weight $\lambda\in P^+$
has the weight space decomposition $V(\lambda)=\sum_{\alpha\in Q^{+}}^{\oplus}
v_{\lambda}U_{\alpha}^+$, where $v_{\lambda}\in V(\lambda)$ is the highest
weight vector of $V(\lambda)$. The weights of the right
$U_q({\mathfrak{g}})$-module $V(\lambda)$ 
coincide with the weights of the left
$U_q({\mathfrak{g}})$-module $V(\lambda)$ and 
the dimensions of the corresponding weight spaces are the same. 

The quantized algebra ${\mathbb{C}}_q\lbrack G \rbrack$ of functions
on the connected simply connected
complex Lie group $G$ with Lie algebra ${\mathfrak{g}}$ is 
the subspace in the linear dual $U_q(\mathfrak{g})^*$ 
spanned by the matrix coefficients 
of the finite-dimensional irreducible representations 
$V(\lambda)$ $(\lambda\in P_+)$.
The Hopf $*$-algebra structure on $U_q(\mathfrak{g})$ induces a Hopf $*$-algebra
structure on ${\mathbb{C}}_q\lbrack G \rbrack\subset U_q({\mathfrak{g}})^*$ by 
the formulas 
\begin{equation}\label{dualstructure}
\begin{split}
&(\phi\psi)(X)=(\phi\otimes\psi)\Delta(X),\quad 1(X)=\varepsilon(X)\\
&\Delta(\phi)(X\otimes Y)=\phi(XY),\quad \varepsilon(\phi)=\phi(1)\\
&S(\phi)(X)=\phi(S(X)),\quad (\phi^*)(X)=\overline{\phi(S(X)^*)},
\end{split}
\end{equation}
where $\phi,\psi\in {\mathbb{C}}_q\lbrack G \rbrack\subset U_q(\mathfrak{g})^*$ 
and $X,Y\in U_q(\mathfrak{g})$.
The algebra ${\mathbb{C}}_q\lbrack G \rbrack$ can be regarded as 
a quantization of the Poisson algebra of polynomial functions
on the algebraic Poisson-Lie group $G$, where the Poisson structure
on $G$ is given by the Sklyanin bracket associated with the 
classical $r$-matrix $-ir$ (cf.\ \eqref{r}).  
Since the $*$-structure \eqref{dualstructure} on ${\mathbb{C}}_q\lbrack G \rbrack$ 
is associated with   
the compact real form $U$ of $G$ in the classical limit, we will write 
${\mathbb{C}}_q\lbrack U \rbrack$ for ${\mathbb{C}}_q\lbrack G \rbrack$ with 
this particular choice
of $*$-structure. 
Note that ${\mathbb{C}}_q\lbrack U \rbrack$ is a $U_q(\mathfrak{g})$-bimodule
with the left respectively right action given by
\begin{equation}\label{lraction}
(X.\phi)(Y):=\phi(YX),\quad (\phi.X)(Y):=\phi(XY) 
\end{equation}
where $\phi\in {\mathbb{C}}_q\lbrack U \rbrack$ and $X,Y\in U_q(\mathfrak{g})$.
The finite-dimensional irreducible 
$U_q(\mathfrak{g})$-module $V(\lambda)$ of highest weight $\lambda\in P_+$
is known to be unitarizable 
(say with inner product $(.,.)$). So we can choose an orthonormal
basis consisting of weight vectors 
\begin{equation}\label{basis}
\lbrace v_{\mu}^{(i)} \, | \, \mu\in P(\lambda), 
i=[1,\hbox{dim}(V(\lambda)_{\mu})]\rbrace,
\end{equation}
where $v_{\mu}^{(i)}\in
V(\lambda)_{\mu}$ (we omit the index $i$ if $\hbox{dim}(V(\lambda)_{\mu})=1$). 
Set
\begin{equation}\label{standardform}
C_{\mu,i;\nu,j}^{\lambda}(X):=(X.v_{\nu}^{(j)},v_{\mu}^{(i)}),\quad X\in 
U_q(\mathfrak{g}),
\end{equation}
for $\mu,\nu\in P(\lambda)$ and $1\leq i\leq\hbox{dim}(V(\lambda)_{\mu})$,
$1\leq j\leq\hbox{dim}(V(\lambda)_{\nu})$. If $\hbox{dim}(V(\lambda)_{\mu})=1$
respectively $\hbox{dim}(V(\lambda)_{\nu})=1$ we omit the dependence on $i$ respectively
$j$ in \eqref{standardform}. It is sometimes also convenient to use
the notation 
\[
C_{v;w}^{\lambda}(X):=(X.w,v), \quad v,w\in V(\lambda), \,\, X\in U_q(\mathfrak{g}).
\]
Note that when $\lambda$ runs through $P_+$ and $\mu,i,\nu$ and $j$ run 
through the above-mentioned sets 
the matrix elements \eqref{standardform} form a linear basis of
${\mathbb{C}}_q\lbrack G \rbrack$. Furthermore, we have the formulas
\begin{equation}\label{relmatrixel}
\begin{split}
\Delta(C_{\mu,i;\nu,j}^{\lambda})&=\sum_{\sigma,s}C_{\mu,i;\sigma,s}^{\lambda}
\otimes C_{\sigma,s;\nu,j}^{\lambda},\\
\varepsilon(C_{\mu,i;\nu,j}^{\lambda})&=\delta_{\mu,\nu}\delta_{i,j},\quad
(C_{\mu,i;\nu,j}^{\lambda})^*=S(C_{\nu,j;\mu,i}^{\lambda}).
\end{split}
\end{equation}
(Sums for which the summation sets are not specified are taken
over the ``obvious'' choice of summation sets). 
Using the relations \eqref{relmatrixel} and
the Hopf algebra axiom for the antipode $S$ we obtain 
\begin{equation}\label{unitarityproperty}
\sum_{\sigma,s}(C_{\sigma,s;\mu,i}^{\lambda})^*C_{\sigma,s;\nu,j}^{\lambda}=
\delta_{\mu,\nu}\delta_{i,j}.
\end{equation}
The elements $(C_{\mu,i;\nu,j}^{\lambda})^*$ are 
matrix coefficients of the contragredient representation $V(\lambda)^*\simeq
V(-\sigma_0\lambda)$ (here $\sigma_0$ is the longest element in $W$).
To be precise, let  $\pi: U_q(\mathfrak{g}) \rightarrow \hbox{End}(V(\lambda))$
be the representation of highest weight $\lambda$, and let 
$(\cdot,\cdot )$ be an inner product with respect to which
$\pi$ is unitarizable.
Fix an orthonormal basis of weight vectors $\lbrace v_{\mu}^{(r)}\rbrace$. 
Let $(\pi^*,V(\lambda)^*)$ be the contragredient
representation, i.e.\  $\pi^*(X)\phi=\phi\circ\pi(S(X))$
for $X\in U_q(\mathfrak{g})$ and $\phi\in V(\lambda)^*$. 
For $u\in V(\lambda)$ set $u^*:=(\cdot,u)\in V(\lambda)^*$.
We define an inner product on $V(\lambda)^*$ by
\[
( u^*,v^*):= \bigl(\pi(K^{-2\rho})v,u\bigr),\quad
u,v\in V(\lambda),
\]
where $\rho=1/2\sum_{\alpha\in R^+}\alpha\in {\mathfrak{h}}^*$. 
By using the fact that $S^2(u)= K^{-2\rho} u K^{2\rho}$ 
($u\in U_q(\mathfrak{g})$) one easily deduces that
$\pi^*$ is unitarizable with respect to the inner product $(\cdot,\cdot)$
on $V(\lambda)^*$ and that 
$\lbrace \phi_{-\mu}^{(i)}:=q^{(\mu,\rho)}(v_{\mu}^{(i)})^*\rbrace$ 
is an orthonormal basis of $V(\lambda)^*$ consisting of weight vectors
(here $\phi_{-\mu}^{(i)}$ has weight $-\mu$).  
Defining the matrix coefficients $C_{-\mu,i;-\nu,j}^{-\sigma_0\lambda}$ of
$(\pi^*,V(\lambda)^*)$ with respect to 
the orthonormal basis $\lbrace \phi_{-\mu}^{(i)}\rbrace$, we then have
\begin{equation}\label{contragred}
(C_{\mu,i;\nu,j}^{\lambda})^*=q^{(\mu-\nu,\rho)}C_{-\mu,i;-\nu,j}^{-\sigma_0\lambda}
\end{equation}
(cf.\ \cite[Prop.\  3.3]{S}).
A fundamental role in Soibel'man's theory of 
irreducible $*$-representations of ${\mathbb{C}}_q\lbrack U \rbrack$
is played by a Poincar{\'e}-Birkhoff-Witt (PBW) type 
factorization of ${\mathbb{C}}_q\lbrack U \rbrack$.
For $\lambda\in P_+$, set
\begin{equation}\label{Blambda}
B_{\lambda}:=\hbox{span}\lbrace C_{v;v_{\lambda}}^{\lambda}
\, | \, v\in V(\lambda)\rbrace.
\end{equation}
Note that $B_{\lambda}$ is a right $U_q(\mathfrak{g})$-submodule
of ${\mathbb{C}}_q\lbrack U \rbrack$ isomorphic to $V(\lambda)$. 
Set
\begin{equation}
A^+:=\bigoplus_{\lambda\in P_+}B_{\lambda}, \quad 
A^{++}:=\bigoplus_{\lambda\in P_{++}} B_{\lambda}.
\end{equation}
The subalgebra and right $U_q({\mathfrak{g}})$-module $A^+$ is equal to
the subalgebra of left $U^+$-invariant elements in
${\mathbb{C}}_q\lbrack U \rbrack$ (cf.\ \cite{J}).
The existence of a PBW type factorization of ${\mathbb{C}}_q\lbrack U \rbrack$ 
now amounts to the following statement.
\begin{Thm}\cite[Theorem 3.1]{S}\label{factorization}
The multiplication map $m: (A^{++})^*\otimes A^{++}
\rightarrow {\mathbb{C}}_q\lbrack U \rbrack$ is surjective.
\end{Thm}
A detailed proof can be found in \cite[Prop.\ 9.2.2]{J}.
The proof is based on certain results concerning
decompositions of tensor products of irreducible 
finite-dimensional $U_q({\mathfrak{g}})$-modules which can be traced back
to Kostant in the classical case \cite[Theorem 5.1]{K}.
The close connection between Theorem \ref{factorization} 
and the decomposition of 
tensor products of irreducible $U_q({\mathfrak{g}})$-modules 
becomes clear by observing that
\begin{equation}\label{tensorrelation}
(B_{\lambda})^*B_{\mu}\simeq V(\lambda)^* \otimes V(\mu)
\end{equation}
as right $U_q(\mathfrak{g})$-modules.

Important for the study of $*$-representations of ${\mathbb{C}}_q\lbrack U \rbrack$
is some detailed information about the
commutation relations between matrix elements in ${\mathbb{C}}_q\lbrack U \rbrack$. 
In view of Theorem \ref{factorization}, 
we are especially interested in commutation relations
between the $C_{\mu,i;\lambda}^{\lambda}$ and $C_{\nu,j;\Lambda}^{\Lambda}$
resp.\ between the $C_{\mu,i;\lambda}^{\lambda}$ and
$(C_{\nu,j;\Lambda}^{\Lambda})^*$, where $\lambda,\Lambda\in P_+$.
To state these commutation relations we need to introduce
certain vector subspaces of $\mathbb{C}_q[U]$.
Let $\lambda,\Lambda\in P_+$ and $\mu\in P(\lambda)$, $\nu\in P(\Lambda)$,
then we set
\begin{equation}\label{N}
\begin{split}
N(\mu,\lambda;\nu,\Lambda)&:=
\hbox{span}\lbrace
C_{v;v_{\lambda}}^{\lambda}C_{w;v_{\Lambda}}^{\Lambda} 
\, | \, (v,w)\in sN\rbrace,\\
N^{opp}(\mu,\lambda;\nu,\Lambda)&:=
\hbox{span}\lbrace
C_{w;v_{\Lambda}}^{\Lambda}C_{v;v_{\lambda}}^{\lambda} 
\, | \, (v,w)\in sN\rbrace
\end{split}
\end{equation}
where $sN:=sN(\mu,\lambda;\nu,\Lambda)$ is the set of pairs
$(v,w)\in V(\lambda)_{\mu'}\times V(\Lambda)_{\nu'}$ with
$\mu'>\mu$, $\nu'<\nu$ and $\mu'+\nu'=\mu+\nu$.
Furthermore set
\begin{equation}\label{O}
\begin{split}
O(\mu,\lambda;\nu,\Lambda)&:=
\hbox{span}\lbrace (C_{v;v_{\lambda}}^{\lambda})^*C_{w;v_{\Lambda}}^{\Lambda}
\, | \, (v,w)\in sO \rbrace,\\
O^{opp}(\mu,\lambda;\nu,\Lambda)&:=
\hbox{span}\lbrace
C_{w;v_{\Lambda}}^{\Lambda}(C_{v;v_{\lambda}}^{\lambda})^*
\, | \, (v,w)\in sO \rbrace
\end{split}
\end{equation} 
where $sO:=sO(\mu,\lambda;\nu,\Lambda)$ is the set of pairs $(v,w)\in
V(\lambda)_{\mu'}\times V(\Lambda)_{\nu'}$ with $\mu'<\mu$, $\nu'<\nu$ and
$\mu-\mu'=\nu-\nu'$.
If $sN$ (resp.\ $sO$) is empty, then let $N=N^{opp}=\lbrace 0 \rbrace$ 
(resp. $O=O^{opp}=\lbrace 0 \rbrace$).
We now have the following proposition.
\begin{Prop}\label{commutation}
Let $\lambda,\Lambda\in P_+$ and $v\in V(\lambda)_{\mu}$, $w\in V(\Lambda)_{\nu}$.\\
{\bf (i)} The matrix elements $C_{v;v_{\lambda}}^{\lambda}$
and $C_{w;v_{\Lambda}}^{\Lambda}$ satisfy the commutation relation
\[
C_{v;v_{\lambda}}^{\lambda}C_{w;v_{\Lambda}}^{\Lambda}=
q^{(\lambda,\Lambda)-(\mu,\nu)}C_{w;v_{\Lambda}}^{\Lambda}
C_{v;v_{\lambda}}^{\lambda}\, \hbox{ mod }\, N(\mu,\lambda;\nu,\Lambda).
\]
Moreover, we have $N=N^{opp}$.\\
{\bf (ii)} The matrix elements  
$(C_{v;v_{\lambda}}^{\lambda})^*$ and $C_{w;v_{\Lambda}}^{\Lambda}$ 
satisfy the commutation relation
\[(C_{v;v_{\lambda}}^{\lambda})^*C_{w;v_{\Lambda}}^{\Lambda}
=q^{(\mu,\nu)-(\lambda,\Lambda)}C_{w;v_{\Lambda}}^{\Lambda}
(C_{v;v_{\lambda}}^{\lambda})^*\, \hbox{ mod }\, O(\mu,\lambda;\nu,\Lambda).\]
Moreover, we have $O=O^{opp}$.
\end{Prop} 
Soibel'man \cite{S} derived commutation relations using the universal $R$-matrix
whereas Joseph \cite[Section 9.1]{J} used the Poincar{\'e}-Birkhoff-Witt 
Theorem for $U_q(\mathfrak{g})$ and the left respectively
right action \eqref{lraction} of 
$U_q(\mathfrak{g})$ on ${\mathbb{C}}_q\lbrack U \rbrack$. 
Although the commutation relations formulated here are slightly sharper,
the proof can be derived in a similar manner and will therefore be omitted.

As a corollary of Proposition \ref{commutation}{\bf (i)} we have
\begin{Cor}\label{other}
Let $\lambda,\Lambda\in P_+$ and $v\in V(\lambda)_{\mu}$, $w\in V(\Lambda)_{\nu}$.
Then
\begin{equation}
C_{v;v_{\lambda}}^{\lambda}C_{w;v_{\Lambda}}^{\Lambda}=
q^{(\mu,\nu)-(\lambda,\Lambda)}C_{w;v_{\Lambda}}^{\Lambda}C_{v;v_{\lambda}}^{\lambda}
\quad \hbox{ mod } \,\, N(\nu,\Lambda;\mu,\lambda).
\end{equation}
\end{Cor}
Note that Proposition \ref{commutation}{\bf (i)} and Corollary \ref{other}
give two  different ways to rewrite 
$C_{v;v_{\lambda}}^{\lambda}C_{w;v_{\Lambda}}^{\Lambda}$ as elements of the 
vector space 
\[
W_{\lambda,\Lambda}:=\hbox{span}\lbrace 
C_{w';v_{\Lambda}}^{\Lambda}C_{v';v_{\lambda}}^{\lambda} 
\, | \, v'\in V(\lambda),\,\, w'\in V(\Lambda) \rbrace.
\]
We will need  both ``inequivalent'' commutation relations (Proposition
\ref{commutation}{\bf (i)} and Corollary \ref{other}) in later sections.
It follows in particular that, when $v'\in V(\lambda)$ and $w'\in V(\Lambda)$ 
run through a basis,
the elements $C_{w';v_{\Lambda}}^{\Lambda}C_{v';v_{\lambda}}^{\lambda}$ 
are (in general) linearly dependent. This also follows  from the following
two observations. On the one hand,
$W_{\lambda,\Lambda}\simeq V(\lambda+\Lambda)$ as right 
$U_q(\mathfrak{g})$-modules. On the other hand,
$V(\lambda+\Lambda)$ occurs with multiplicity one in
$V(\lambda)\otimes V(\Lambda)$,
whereas in general $V(\lambda)\otimes V(\Lambda)$
has other irreducible components too.

By contrast, the commutation relation given in Proposition
\ref{commutation}{\bf (ii)} is unique in the sense that,
when $v\in V(\lambda)$ and $w\in V(\Lambda)$ run through a basis,
the $C_{w;v_{\Lambda}}^{\Lambda}(C_{v;v_{\lambda}}^{\lambda})^*$ 
are linearly independent (cf.\ \eqref{tensorrelation}).

We end this section by recalling the special case
$\mathfrak{g}=\mathfrak{s}\mathfrak{l}(2,{\mathbb{C}})$. Set
\begin{equation}
\begin{split}
&t_{11}:=C_{\varpi_1;\varpi_1}^{\varpi_1},\quad
 t_{12}:=C_{\varpi_1;-\varpi_1}^{\varpi_1},\\
&t_{21}:=C_{-\varpi_1;\varpi_1}^{\varpi_1},\quad
 t_{22}:=C_{-\varpi_1;-\varpi_1}^{\varpi_1}.
\end{split}
\end{equation} 
Then it is well known that the $t_{ij}$'s generate the algebra
${\mathbb{C}}_q\lbrack SU(2) \rbrack$. The commutation relations
\begin{equation}
\begin{split}
&t_{k1}t_{k2}=qt_{k2}t_{k1},\quad t_{1k}t_{2k}=qt_{2k}t_{1k} \quad (k=1,2),\\
&t_{12}t_{21}=t_{21}t_{12},\quad t_{11}t_{22}-t_{22}t_{11}=
(q-q^{-1})t_{12}t_{21},\\
&t_{11}t_{22}-qt_{12}t_{21}=1
\end{split}
\end{equation}
characterize the algebra structure of 
$\mathbb{C}_q[SU(2)]$ in terms of the generators
$t_{ij}$. The $*$-structure is uniquely determined by the formulas
$t_{11}^*=t_{22}$, $t_{12}^*=-qt_{21}$.
\section{Quantized function algebras on generalized flag manifolds}   
\label{section:quan}
Let $S$ be any subset of the simple roots $\Delta$. 
We will sometimes identify $S$ with the index set 
$\lbrace i \, | \, \alpha_i\in S\rbrace$.
Let $\mathfrak{p}_S\subset \mathfrak{g}$ be the corresponding
standard parabolic subalgebra, given explicitly by \eqref{pS}.
We define the quantized universal enveloping algebra $U_q({\mathfrak{l}}_S)$
associated with the Levi factor ${\mathfrak{l}}_S$ of ${\mathfrak{p}}_S$
as the subalgebra of $U_q({\mathfrak{g}})$ generated
by $K_i^{\pm 1}$ ($i\in [1,r]$) and $X_i^{\pm}$ ($i\in S$). Note that
$U_q({\mathfrak{l}}_S)$ is a Hopf $*$-subalgebra of $U_q({\mathfrak{g}})$.  

For later use in this section we briefly discuss the finite-dimensional
representation theory of $U_q({\mathfrak{l}}_S)$. Recall that 
${\mathfrak{l}}_S$ is a reductive Lie algebra with centre
\begin{equation}
Z({\mathfrak{l}}_S)=\bigcap_{i\in S}\hbox{Ker}(\alpha_i)\subset {\mathfrak{h}}.
\end{equation}
Moreover, we have direct sum decompositions 
\begin{equation}\label{decompCartan}
{\mathfrak{h}}=Z({\mathfrak{l}}_S)\oplus {\mathfrak{h}}_S,\quad
{\mathfrak{l}}_S=Z({\mathfrak{l}}_S)\oplus {\mathfrak{l}}_S^0,
\end{equation}
where ${\mathfrak{h}}_S=\hbox{span}\lbrace H_i \rbrace_{i\in S}$ and
${\mathfrak{l}}_S^0$ is the semisimple part of ${\mathfrak{l}}_S$. 
The semisimple part ${\mathfrak{l}}_S^0$ is explictly given by
\begin{equation}
{\mathfrak{l}}_S^0:={\mathfrak{h}}_S\oplus\bigoplus_{\alpha\in \Gamma_S\cap
(-\Gamma_S)}
{\mathfrak{g}}_{\alpha}.
\end{equation}
We define the quantized universal enveloping algebra $U_q({\mathfrak{l}}_S^0)$
associated with the semisimple part ${\mathfrak{l}}_S^0$ of ${\mathfrak{l}}_S$
as the subalgebra of $U_q({\mathfrak{g}})$ generated
by $K_i^{\pm 1}$ and  $X_i^{\pm}$  for all $i\in S$. Observe that 
$U_q({\mathfrak{l}}_S^0)$ is a Hopf $*$-subalgebra of $U_q({\mathfrak{g}})$.
\begin{Prop}
Any finite-dimensional $U_q({\mathfrak{l}}_S)$-module $V$ which is completely
reducible as $U_q({\mathfrak{h}})$-module, is completely reducible as
$U_q({\mathfrak{l}}_S)$-module. 
\end{Prop}
\begin{proof}
Let $V$ be a finite-dimensional 
left $U_q({\mathfrak{l}}_S)$-module which is completely
reducible as $U_q({\mathfrak{h}})$-module. Then the linear subspace
\[V^+:=\lbrace v\in V \, | \, X_i^+v=0 \quad \forall i\in S \rbrace\]
is $U_q({\mathfrak{h}})$-stable and splits as a direct sum of weight spaces.
Let $\{ v_i\}$ be a linear basis of $V^+$ consisting of weight vectors, and 
set $V_i:=U_q({\mathfrak{l}}_S^0)v_i$. Since $U_q({\mathfrak{l}}_S^0)$ is the quantized
universal enveloping algebra associated with a semisimple Lie algebra, it follows
that $V=\sum_i^{\oplus}V_i$ is a decomposition of $V$ into irreducible 
$U_q({\mathfrak{l}}_S^0)$-modules. On the other hand, the 
$V_i$ are $U_q({\mathfrak{l}}_S)$-stable
since the vectors $v_i$ are weight vectors. Hence $V=\sum_i^{\oplus}V_i$
is a decomposition of $V$ into irreducible 
$U_q({\mathfrak{l}}_S)$-modules.
\end{proof} 
There are obvious notions of weight vectors and weights for 
$U_q({\mathfrak{l}}_S)$-modules.
With a suitably extended interpretation of the notion of highest weight,
the irreducible finite-dimensional $U_q({\mathfrak{l}}_S)$-modules
may be characterized in terms of highest weights. 
We shall only be interested in irreducible
$U_q({\mathfrak{l}}_S)$-modules with weights in the lattice $P$.
For instance, the restriction of an irreducible $P$-weighted 
$U_q(\mathfrak{g})$-module to $U_q(\mathfrak{l}_S)$ decomposes 
into such irreducible $U_q(\mathfrak{l}_S)$-modules.

Branching rules for the restriction of finite-dimensional representations of
$U_q(\mathfrak{g})$ to $U_q(\mathfrak{l}_S)$ are determined by the
behaviour of the corresponding characters. Since the characters for
$P$-weighted irreducible finite-dimensional
representations of $U_q(\mathfrak{g})$ and $U_q(\mathfrak{l}_S)$
are the same as for the corresponding representations of $\mathfrak{g}$ and
$\mathfrak{l}_S$, one easily derives the following proposition.
\begin{Prop}\label{zelfde}
Let $\lambda\in P_+$. The multiplicity of any $P$-weighted irreducible
$U_q(\mathfrak{l}_S)$-module in the irreducible decomposition of the restriction of 
the $U_q(\mathfrak{g})$-module $V(\lambda)$ to 
$U_q(\mathfrak{l}_S)$ is the same as in the classical case.
\end{Prop}
Next, we define the quantized algebra of functions on $U/K_S$.
The mapping $\iota_S^*\colon U_q({\mathfrak{g}})^*\twoheadrightarrow
U_q({\mathfrak{l}}_S)^*$ dual to the Hopf $*$-embedding 
$\iota_S: U_q({\mathfrak{l}}_S)\hookrightarrow U_q({\mathfrak{g}})$
is surjective, and we set
\[
{\mathbb{C}}_q\lbrack L_S\rbrack:=\iota_S^*({\mathbb{C}}_q\lbrack G \rbrack)=
\lbrace\phi\circ \iota_S \, | \, \phi\in {\mathbb{C}}_q\lbrack G \rbrack\rbrace.
\]
The formulas \eqref{dualstructure} uniquely determine a Hopf $*$-algebra structure on
${\mathbb{C}}_q\lbrack L_S\rbrack$, and $\iota_S^*$ then becomes a  Hopf $*$-algebra
morphism. We write ${\mathbb{C}}_q[K_S]$ for ${\mathbb{C}}_q[L_S]$ 
with this particular choice of $*$-structure.
Assume now that $S\not=\Delta$.  
Define a $*$-subalgebra ${\mathbb{C}}_q\lbrack U/K_S \rbrack \subset
{\mathbb{C}}_q\lbrack U \rbrack$ by
\begin{equation}
\begin{split}
{\mathbb{C}}_q\lbrack U/K_S \rbrack&:=\lbrace \phi\in {\mathbb{C}}_q\lbrack U \rbrack 
\, | \, (\hbox{id}\otimes \iota_S^*)\Delta(\phi)=\phi\otimes 1\rbrace\\
&=\lbrace \phi\in {\mathbb{C}}_q\lbrack U \rbrack \, | \, X.\phi=\varepsilon(X)\phi,
\quad \forall\, X\in U_q({\mathfrak{l}}_S)\rbrace
\end{split}
\end{equation}
The algebra ${\mathbb{C}}_q\lbrack U/K_S \rbrack$ is a 
left ${\mathbb{C}}_q\lbrack U \rbrack$-subcomodule of
${\mathbb{C}}_q\lbrack U \rbrack$. We call it
the quantized algebra of functions on the generalized flag manifold $U/K_S$.

In a similar way, one can define the quantized function algebra
${\mathbb{C}}_q\lbrack K^0_S \rbrack$ corresponding to the 
semisimple part $K_S^0$ of $K_S$ as the image of the dual of the
natural embedding $U_q({\mathfrak{l}}_S^0)\hookrightarrow U_q({\mathfrak{g}})$.
Its Hopf $*$-algebra structure is again given by the formulas \eqref{dualstructure}.
The subalgebra ${\mathbb{C}}_q\lbrack U/K_S^0\rbrack$ then consists by definition
of all right ${\mathbb{C}}_q\lbrack K_S^0\rbrack$-invariant elements
in ${\mathbb{C}}_q\lbrack U\rbrack$.
Note that ${\mathbb{C}}_q\lbrack U/K_S^0\rbrack\subset {\mathbb{C}}_q\lbrack U\rbrack$ 
is a left $U_q({\mathfrak{h}})$-submodule
and that ${\mathbb{C}}_q\lbrack U/K_S\rbrack$ 
coincides with the subalgebra of $U_q({\mathfrak{h}})$-invariant elements
in ${\mathbb{C}}_q\lbrack U/K_S^0\rbrack$.

We now turn to PBW type factorizations of the algebra
${\mathbb{C}}_q\lbrack U/K_S\rbrack$.
Write $P(S)$, $P_+(S)$, resp.\ $P_{++}(S)$ for
${\mathbb{K}}\hbox{-span}\lbrace\varpi_{\alpha}\rbrace_{\alpha\in S}$
with ${\mathbb{K}}={\mathbb{Z}}$, ${\mathbb{Z}}_+$ resp.\  ${\mathbb{N}}$.
Set $S^c:=\Delta\setminus S$.
The quantized algebra $A_S^{hol}$ of holomorphic polynomials on $U/K_S$ is defined by
\begin{equation}\label{AShol}
A_S^{hol}:=\bigoplus_{\lambda\in P_+(S^c)}B_{\lambda}\subset {\mathbb{C}}_q\lbrack
U\rbrack,
\end{equation}
where $B_{\lambda}$ is given by \eqref{Blambda} (cf.\ \cite{LR1}, \cite{LR2},
\cite{S2}, \cite{Jur} and \cite{Kor}). 
Note that $A_S^{hol}$ is a right $U_q({\mathfrak{g}})$-comodule subalgebra of 
${\mathbb{C}}_q\lbrack U \rbrack$, \eqref{AShol} being the
(multiplicity free) decomposition of $A_S^{hol}$ into irreducible
$U_q({\mathfrak{g}})$-modules. 
The right $U_q({\mathfrak{g}})$-module algebra $(A_S^{hol})^*\subset 
{\mathbb{C}}_q\lbrack U\rbrack$ is called the quantized
algebra of antiholomorphic polynomials on $U/K_S$.
\begin{Lem}
The linear subspace   
\[
A_S^0:=m\bigl((A_S^{hol})^*\otimes A_S^{hol}\bigr)\subset {\mathbb{C}}_q\lbrack
U\rbrack,
\]
where $m$ is the multiplication map of ${\mathbb{C}}_q\lbrack U \rbrack$, is 
a right $U_q({\mathfrak{g}})$-submodule $*$-subalgebra of 
${\mathbb{C}}_q\lbrack U\rbrack$. 
\end{Lem}
\begin{proof}
Proposition \ref{commutation}{\bf (ii)} implies that 
$A_S^0$ is a subalgebra of ${\mathbb{C}}_q\lbrack U \rbrack$. 
The other assertions are immediate.
\end{proof}
The subalgebra $A_S^0$ may be considered as a quantum analogue of
the algebra of complex-valued polynomial functions on the real manifold $U/K_S^0$.
\begin{rem}\label{Litremark} 
In the classical setting ($q=1$), 
the algebra $A_S^0$ ($\#S^c=1$) can be interpreted 
as algebra of functions on the product of an affine
spherical $G$-variety with its dual. The $G$-module structure on $A_S^0$ is then 
related to the doubled $G$-action (see \cite{Pan1}, \cite{Pan} for the terminology).  
These (and related) $G$-varieties have been studied in several papers, see for example
\cite{Pan}, \cite{Pan1} and \cite{Lit1}.
\end{rem}
The algebra $A_S^0\subset {\mathbb{C}}_q\lbrack U\rbrack$ is stable under the left
$U_q({\mathfrak{h}})$-action, so we can speak of $U_q({\mathfrak{h}})$-weighted 
elements in $A_S^0$.
Let $A_S$ be the left $U_q({\mathfrak{h}})$-invariant elements of $A_S^0$.
Then $A_S\subset {\mathbb{C}}_q\lbrack U\rbrack$ 
is a right $U_q({\mathfrak{g}})$-module $*$-subalgebra of
${\mathbb{C}}_q\lbrack U\rbrack$. 
We now have the following lemma.
\begin{Lem}\label{definitionAS}
We have $A_S^0\subset {\mathbb{C}}_q\lbrack U/K_S^0\rbrack$, 
so in particular $A_S\subset {\mathbb{C}}_q\lbrack U/K_S\rbrack$.
Furthermore,
\begin{equation}\label{matrixvorm}
A_S=\hbox{span}\lbrace (C_{v;v_{\lambda}}^{\lambda})^*C_{w;v_{\lambda}}^{\lambda}
\, | \, \lambda\in P_+(S^c),\,\, v,w\in V(\lambda) \rbrace.
\end{equation}
\end{Lem}
\begin{proof}
Choose $\lambda\in P_+(S^c)$ and $i\in S$. Then we have $X_i^+\cdot v_{\lambda}=0$ and
$K_i\cdot v_{\lambda}=v_{\lambda}$. It follows that 
$\mathbb{C} v_\lambda \subset V(\lambda)$ is a one-dimensional 
$U_{q_i}(\mathfrak{s}\mathfrak{l}(2;{\mathbb{C}}))$-submodule, 
where we consider the $U_{q_i}(\mathfrak{s}\mathfrak{l}(2;{\mathbb{C}}))$ action 
on $V(\lambda)$ via the embedding 
$\phi_i\colon U_{q_i}(\mathfrak{s}\mathfrak{l}(2;{\mathbb{C}}))
\hookrightarrow U_q({\mathfrak{g}})$. 
It follows that $X_i^{-}\cdot v_{\lambda}=0$.
This readily implies that $A_S^0\subset {\mathbb{C}}_q\lbrack U/K_S^0 \rbrack$. 
The remaining assertions are immediate.
\end{proof}
\begin{Def}
We call $A_S\subset {\mathbb{C}}_q\lbrack U/K_S \rbrack$ the factorized 
$*$-subalgebra associated with $U/K_S$.
\end{Def}
In view of Theorem \ref{factorization}, there is reason to expect 
that the factorized algebra $A_S$ is equal to ${\mathbb{C}}_q\lbrack
U/K_S\rbrack$ for any generalized flag manifold $U/K_S$. 
Although we cannot prove this in general, we do have a proof (cf.\ Theorem \ref{yes}) 
for a certain subclass of generalized flag manifolds that we shall define and classify
in the following proposition. 
For the proof in these cases we use the so-called Parthasarathy-Ranga
Rao-Varadarajan (PRV) conjecture, which was proved independently 
by Kumar \cite{Ku} and Mathieu \cite{Ma}. The PRV conjecture gives information 
about which irreducible constituents occur in tensor products of irreducible 
finite-dimensional ${\mathfrak{g}}$-modules.
It seems likely that a proof for arbitrary generalized flag manifold $U/K_S$
would require further detailed information about irreducible decompositions of  
tensor products of finite-dimensional representations of $\mathfrak{g}$.

Recall the notations introduced in section 2.
The following proposition was observed by Koornwinder \cite{Koo}.
\begin{Prop}\label{Koornwinderpairs}
\textup{(} \cite{Koo}\textup{)}
Let $U$ be a connected, simply connected compact Lie group with Lie algebra
${\mathfrak{u}}$,  and let ${\mathfrak{p}}\subset {\mathfrak{g}}$ be 
a standard maximal parabolic subalgebra. 
Let $K\subset U$ be the connected subgroup with Lie algebra
${\mathfrak{k}}:={\mathfrak{p}}\cap {\mathfrak{u}}$. 
Then $(U,K)$ is a Gel'fand pair if and only if one of the 
following three conditions are satisfied:\\
{\bf (i)} $(U,K)$ is an irreducible compact Hermitian symmetric pair;\\
{\bf (ii)} $(U,K)\simeq  (SO(2l+1),U(l))$, \quad $(l\geq 2)$;\\
{\bf (iii)} $(U,K)\simeq (Sp(l),U(1)\times Sp(l-1))$, 
\quad $(l\geq 2)$.
\end{Prop}
\begin{proof}
For a list of the irreducible compact Hermitian symmetric pairs  
see \cite[Ch. X, Table V]{H}. The proposition 
follows from this and the classification of compact Gelfand pairs $(U,K)$ with 
$U$ simple (cf.\ \cite[Tabelle 1]{Kr}). 
\end{proof}
Let $(U,K)$ be a pair  
from the list {\bf (i)}--{\bf (iii)} in Proposition \ref{Koornwinderpairs}, and let
$({\mathfrak{u}}, {\mathfrak{k}})$ be the associated pair of Lie algebras. 
Then  ${\mathfrak{k}}={\mathfrak{k}}_S$ for some subset $S\subset \Delta$
with $\#S^c=1$. We call the simple root $\alpha\in S^c$ the Gel'fand node associated
with $(U,K)$. 
\begin{Prop}\label{Koor}
Let $(U,K)$ be a  pair  
from the list {\bf (i)}--{\bf (iii)} in Proposition \ref{Koornwinderpairs},
and let $\alpha\in \Delta$ be the associated Gel'fand node with
corresponding fundamental weight $\varpi:=\varpi_{\alpha}$.
Let $\lbrace \mu_1,\ldots,\mu_l\rbrace$
be the fundamental spherical weights of $(U,K)$.
Then every fundamental spherical representation $V(\mu_i)$ occurs 
in the decomposition of $V(\varpi)^*\otimes V(\varpi)$.
\end{Prop}
\begin{proof}
For the proof we use the PRV conjecture, which states the following.
Let $\lambda,\mu\in P_+$ and $w\in W$.
Let $[\lambda+w\mu]$ be the unique element in $P_+$ which lies in the
$W$-orbit of $\lambda+w\mu$. Then $V([\lambda+w\mu])$ occurs with
multiplicity at least one in $V(\lambda)\otimes V(\mu)$.
The procedure is now a follows. For a pair $(U,K)$ from the list {\bf (i)}--{\bf
(iii)} of Proposition \ref{Koornwinderpairs} 
we write down the fundamental spherical weights $\lbrace \mu_i\rbrace_{i=1}^l$
in terms of the fundamental weights $\lbrace \varpi_i\rbrace_{i=1}^r$
(cf.\  \cite[Tabelle 1]{Kr}, or in the case of Hermitian symmetric spaces
one can also write them down from the corresponding Satake diagrams \cite{Su}).
Then we look for Weyl group elements $w_i\in W$ such that
\[
[\varpi-w_i\varpi]=\mu_i,\quad (i=[1,l])
\]
(here we used that $V(\varpi)^*\simeq V(-\sigma_0\varpi)$).

As an example, let us follow the procedure for the compact Hermitian symmetric
pair $(U,K)=(SO(2l),U(l))$ ($l\geq 2$).
We use the standard realization of the root system $R$ of type $D_l$ 
in the $l$-dimensional vector space $V=\sum_{i=1}^l{\mathbb{R}}\varepsilon_i$, with
basis given by $\alpha_i=\varepsilon_i-\varepsilon_{i+1}$ ($i=[1,l-1]$)
and $\alpha_l=\varepsilon_{l-1}+\varepsilon_l$.  
The fundamental weights are given by
\begin{equation}
\begin{split}
\varpi_i&=\varepsilon_1+\varepsilon_2+\ldots +\varepsilon_i,\quad (i<l-1),\nonumber\\
\varpi_{l-1}&=(\varepsilon_1+\varepsilon_2+\ldots +
\varepsilon_{l-1}-\varepsilon_l)/2,\nonumber\\
\varpi_l&=(\varepsilon_1+\varepsilon_2+\ldots + 
\varepsilon_{l-1}+\varepsilon_l)/2.\nonumber
\end{split}
\end{equation}
We set ${\varpi}=\varpi_l$ (i.e. $S^c=\lbrace \alpha_l\rbrace$).
Let $\sigma_i$ be the linear map defined by $\varepsilon_j\mapsto -\varepsilon_j$ 
($j=i,i+1$) and $\varepsilon_j\mapsto\varepsilon_j$ otherwise.
Then $\sigma_i\in W$ ($i=[1,l-1]$). If $l=2l'+1$, then
\begin{equation}\label{fundweightodd}
\begin{split}
\varpi-\sigma_1\sigma_3\ldots\sigma_{2i-1}\varpi
&=\varpi_{2i},\quad (i=[1,l'-1]),\\ 
\varpi-\sigma_1\sigma_3\ldots\sigma_{2l'-1}\varpi
&=\varpi_{l-1}+\varpi_l.
\end{split}
\end{equation}
If $l=2l'$ then we have
\begin{equation}\label{fundweighteven}
\begin{split}
\varpi-\sigma_1\sigma_3\ldots\sigma_{2i-1}\varpi
&=\varpi_{2i},\qquad (i=[1,l'-1]),\\
\varpi-\sigma_1\sigma_3\ldots\sigma_{2l'-1}\varpi
&=2\varpi_l.
\end{split}
\end{equation}
By comparison with \cite[Tabelle 1]{Kr} 
we see from  \eqref{fundweightodd} (resp.\ \eqref{fundweighteven})
that all the fundamental spherical weights of the pair
$(U,K)=(SO(2l),U(l))$ have been obtained. 
The other cases are checked in a similar manner.
\end{proof}
The question naturally arises whether the fundamental spherical
representations occur with multiplicity one in $V(\varpi)^*\otimes V(\varpi)$
and whether they exhaust (together with the trivial representation) 
the irreducible components of $V(\varpi)^*\otimes V(\varpi)$. 
For the complex Grassmannians 
$(U,K)=\bigl(SU(p+l),S(U(p)\times U(l))\bigr)$ this is
indeed the case (this can be easily proved using the Pieri formula 
for Schur functions \cite[Chapter I, (5.17)]{M}, see \cite{DS} for more details).
For the general case we do not known the answer, but it is true that
the irreducible decomposition of $V(\varpi)^*\otimes V(\varpi)$
is multiplicity free and that all irreducible components
are spherical representations. This is an easy consequence of the
proof of the following main theorem of this section.
\begin{Thm}\label{yes}
The factorized $*$-subalgebra $A_S$ is equal to ${\mathbb{C}}_q\lbrack U/K_S \rbrack$ 
if\\
{\bf (i)} $S=\emptyset$, i.e. $U/K_S=U/T$ is the full flag manifold;\\
{\bf (ii)} $\#S^c=1$ and the simple root $\alpha\in S^c$ is a Gel'fand node.
\end{Thm}
\begin{proof}
To prove {\bf (i)} we look at the simultaneous eigenspace decomposition of 
${\mathbb{C}}_q\lbrack U \rbrack$ with 
respect to the left $U_q({\mathfrak{h}})$-action on ${\mathbb{C}}_q\lbrack U \rbrack$.
The simultaneous eigenspace corresponding to the character
$\varepsilon$ on $U_q({\mathfrak{h}})$ is exactly ${\mathbb{C}}_q\lbrack U/T\rbrack$. 
Using Soibel'man's 
factorization of ${\mathbb{C}}_q\lbrack U \rbrack$ 
(cf.\ Theorem \ref{factorization}) and Lemma \ref{definitionAS}, 
it is then easily checked 
that ${\mathbb{C}}_q\lbrack U/T\rbrack=A_{\emptyset}$. 
To prove {\bf (ii)} we note that 
\[
\bigoplus_{i=1}^l V(\mu_i)\hookrightarrow V(\varpi)^*\otimes
V(\varpi)\simeq (B_{\varpi})^*B_{\varpi}\subset A_S 
\]
as right $U_q({\mathfrak{g}})$-modules by Proposition \ref{Koor}
(here we use the notations as introduced in Proposition \ref{Koor}). 
Now  ${\mathbb{C}}_q\lbrack U \rbrack$ is an integral domain
(cf.\  \cite[Lemma 9.1.9 {\bf (i)}]{J}), hence
$v_{\lambda}v_{\mu}\in A_S$
is a highest weight vector of highest weight $\lambda+\mu$ if $v_{\lambda},v_{\mu}\in
A_S$ are highest weight vectors of highest weight $\lambda$ respectively $\mu$.
It follows that
\[\bigoplus_{n_i\in {\mathbb{Z}}_{+}}
V(n_1\mu_1+\ldots n_l\mu_l)\hookrightarrow A_S\]
as right $U_q({\mathfrak{g}})$-modules.
On the other hand we have the decomposition 
\[{\mathbb{C}}_q\lbrack U/K_S \rbrack\simeq \bigoplus_{n_i\in {\mathbb{Z}}_{+}}
V(n_1\mu_1+\ldots n_l\mu_l)
\]
of ${\mathbb{C}}_q\lbrack U/K_S \rbrack$ into irreducible
right $U_q({\mathfrak{g}})$-modules by
Proposition \ref{zelfde}. This implies $A_S={\mathbb{C}}_q\lbrack U/K_S \rbrack$. 
\end{proof}
In the remainder of the paper we study the irreducible $*$-representations 
of the $*$-algebras $A_S$ and ${\mathbb{C}}_q[U/K_S]$. In the next section we first
consider the restriction of the irreducible $*$-representations of ${\mathbb{C}}_q[U]$
to the $*$-algebras $A_S$ and ${\mathbb{C}}_q[U/K_S]$.
\section{Restriction of irreducible $*$-representations 
to ${\mathbb{C}}_q\lbrack U/K \rbrack$}
Let us first recall some results from
Soibel'man \cite{S} concerning the irreducible
$*$-representations of ${\mathbb{C}}_q\lbrack U\rbrack$.
Let $\lbrace e_i\rbrace_{i\in {\mathbb{Z}}_+}$ 
be the standard orthonormal basis of $l_2({\mathbb{Z}}_+)$.  
Write ${\mathcal{B}}(l_2({\mathbb{Z}}_+))$ for the algebra of
bounded linear operators on $l_2({\mathbb{Z}}_+)$.
Then the formulas
\begin{equation}\label{action2}
\begin{split}
&\pi_q(t_{11})e_j=\sqrt{(1-q^{2j})}e_{j-1},\quad \pi_q(t_{12})e_j=-q^{j+1}e_j,\\
&\pi_q(t_{21})e_j=q^je_j,\quad \pi_q(t_{22})e_j=\sqrt{(1-q^{2(j+1)})}e_{j+1}\\
\end{split}
\end{equation}
(here $\pi_q(t_{11})e_0=0$)
uniquely determine an irreducible $*$-representation 
\[\pi_q: {\mathbb{C}}_q\lbrack SU(2) \rbrack\rightarrow
{\mathcal{B}}(l_2({\mathbb{Z}}_+)).\]
Now the dual of the injective Hopf $*$-algebra morphism $\phi_i\colon 
U_{q_i}(\mathfrak{s}\mathfrak{l}(2;{\mathbb{C}}))\hookrightarrow 
U_q(\mathfrak{g})$ corresponding to the $i$th node
of the Dynkin diagram ($i\in [1,r]$) is a surjective Hopf $*$-algebra morphism
$\phi_i^*\colon \mathbb{C}_q\lbrack U \rbrack\twoheadrightarrow 
\mathbb{C}_{q_i}\lbrack SU(2) \rbrack$. Hence we obtain
irreducible $*$-representations 
$\pi_i:=\pi_{q_i}\circ\phi_i^*:{\mathbb{C}}_q\lbrack U \rbrack\rightarrow
{\mathcal{B}}(l_2({\mathbb{Z}}_+))$. 

On the other hand, there is a family of one-dimensional 
$*$-representations $\tau_t$ of ${\mathbb{C}}_q\lbrack U \rbrack$
parametrized by the maximal torus $t\in T\simeq \mathbb{T}^r$ 
($\mathbb{T}\subset \mathbb{C}$
denoting the unit circle in the complex plane). More explicitly, let
$\iota_T:U_q(\mathfrak{h})\hookrightarrow U_q(\mathfrak{g})$ be the natural
Hopf $*$-algebra embedding, and set 
${\mathbb{C}}_q\lbrack T\rbrack:=
\hbox{span}\lbrace \phi_{\mu}\rbrace_{\mu\in P}\subset U_q(\mathfrak{h})^*$,
where $\phi_{\mu}(K^{\sigma}):=q^{(\mu,\sigma)}$ for $\sigma\in Q$. 
As in \eqref{dualstructure} we get a 
Hopf $*$-algebra structure on ${\mathbb{C}}_q\lbrack T\rbrack$.
Then $\iota_T^*: {\mathbb{C}}_q\lbrack U \rbrack\rightarrow 
{\mathbb{C}}_q\lbrack T\rbrack$, $\iota_T^*(\phi):=\phi\circ\iota_T$
is a surjective Hopf $*$-algebra morphism. 
Any irreducible $*$-representation of ${\mathbb{C}}_q\lbrack T\rbrack$ is 
one-dimensional and can be written 
as $\tilde{\tau}_t(\phi_{\mu}):=t^{\mu}$ 
for a unique $t\in T\simeq \mathbb{T}^r$.
Here $t^{\mu}:=t_1^{m_1}\ldots t_r^{m_r}$ for
$\mu=\sum_{i=1}^rm_i\varpi_i$. 
So we obtain a one-dimensional
$*$-representation $\tau_t:=\tilde{\tau}_t\circ\iota_T^*$ 
of ${\mathbb{C}}_q\lbrack U \rbrack$, which is given explicitly 
on matrix elements $C_{\mu,i;\nu,j}^{\lambda}$ by the formula
\begin{equation}\label{one-dimensional}
\tau_t(C_{\mu,i;\nu,j}^{\lambda})=\delta_{\mu,\nu}\delta_{i,j}t^{\mu}.
\end{equation}
The following theorem completely describes the irreducible $*$-representations
of ${\mathbb{C}}_q\lbrack U \rbrack$.
\begin{Thm}[Soibel'man \cite{S}]\label{Soibelman}
Let $\sigma\in W$, and fix a reduced 
expression $\sigma=s_{i_1}s_{i_1}\cdots s_{i_l}$.
The $*$-representation 
\begin{equation}
\pi_\sigma:=\pi_{i_1}\otimes\pi_{i_2}\otimes\cdots\otimes\pi_{i_l}
\end{equation}
does not depend on the choice of reduced expression 
\textup{(}up to equivalence\textup{)}. The set
\[
\lbrace \pi_\sigma\otimes\tau_t \, | \, t\in T, \sigma\in W \rbrace
\]
is a complete set of mutually inequivalent irreducible $*$-representations of
${\mathbb{C}}_q\lbrack U \rbrack$.
\end{Thm}
Here tensor products of $*$-representations are defined in the usual way
by means of the coalgebra structure on ${\mathbb{C}}_q\lbrack U \rbrack$.
The irreducible representation $\pi_e$ with respect to the unit element $e\in W$
is the one-dimensional $*$-representation associated with the counit $\epsilon$ on
${\mathbb{C}}_q[U]$.
In Soibel'man's terminology, the representations $\pi_{\sigma}\otimes \tau_t$
are said to be associated with the Schubert cell $X_{\sigma}$ of $U/T$ 
(cf.\ section 2).

We also mention here an important property of the kernel of $\pi_{\sigma}$, 
which we will repeatedly need later on. Let $U_q({\mathfrak{b}}_+)$ be 
the subalgebra of $U_q({\mathfrak{g}})$ generated by the $K_i^{\pm 1}$ 
and the $X_i^{+}$ ($i\in [1,r]$). For any $\lambda\in P_+$,
the $*$-representation $\pi_\sigma$ satisfies
\begin{equation}\label{part1}
\pi_\sigma(C_{v;v_{\lambda}}^{\lambda})=0\quad \bigl(v\notin 
U_q({\mathfrak{b}}_+)v_{\sigma\lambda}\bigr),\quad \pi_{\sigma}(C_{v_{\sigma\lambda};
v_{\lambda}}^{\lambda})\not=0
\end{equation}
(cf.\  \cite[Theorem 5.7]{S}). Formula \eqref{part1} combined with 
\cite[Lemma 2.12]{BGG} shows that the classical limit of
the kernel of $\pi_\sigma$ formally tends to the ideal of 
functions vanishing on $X_\sigma$.  

Fix now a subset $S\subsetneq \Delta$. We freely use the notations 
introduced earlier. Our next goal is to 
describe how the $\ast$-representations $\pi_\sigma$ decompose
under restriction to the subalgebra ${\mathbb{C}}_q\lbrack U/K_S \rbrack$.
Consider the selfadjoint operators
\begin{equation}\label{Loperators}
L_{\sigma\lambda;\lambda}:=
\pi_\sigma((C_{\sigma\lambda;\lambda}^{\lambda})^*C_{\sigma\lambda;\lambda}^{\lambda})
\end{equation}
for $\lambda\in P_{+}(S^c)$. 
Let $\sigma=s_{i_1}\cdots s_{i_l}$ be a reduced expression for $\sigma$, and set
$\pi_\sigma=\pi_{i_1}\otimes \pi_{i_2}\otimes\cdots\otimes \pi_{i_l}$.
Then it follows from  \cite[Proof of Prop.\  5.2]{S} (see also \cite[Proof of
Prop.\ 5.8]{S}) that
\begin{equation}\label{facto}
\pi_\sigma(C_{\sigma\lambda;\lambda}^{\lambda})=
c\,\pi_{q_{i_1}}(t_{21})^{(\lambda,\gamma_1\spcheck)}\otimes 
\pi_{q_{i_2}}(t_{21})^{(\lambda,\gamma_2\spcheck)}\otimes\cdots\otimes
\pi_{q_{i_l}}(t_{21})^{(\lambda,\gamma_l\spcheck)}
\end{equation}
where the scalar $c\in {\mathbb{T}}$ depends on the particular choices
of bases for the irreducible representations $V(\mu)$ ($\mu\in P_+$), and with
\begin{equation}\label{factohelp}
\gamma_k:=s_{i_l}s_{i_{l-1}}\cdots s_{i_{k+1}}(\alpha_{i_k})\quad (1\leq k\leq l-1),
\quad \gamma_l:=\alpha_{i_l}.
\end{equation}
The proof of \eqref{Loperators}, which was given in \cite{S} under the assumption
that $\lambda\in P_{++}$, is in fact valid for all dominant weights $\lambda\in P_+$.
It follows from  \eqref{action2}, \eqref{Loperators} and \eqref{facto} 
that $l_2({\mathbb{Z}}_+)^{\otimes l(\sigma)}$ 
decomposes as an orthogonal direct sum of eigenspaces for $L_{\sigma\lambda;\lambda}$,
\begin{equation}
l_2({\mathbb{Z}}_+)^{\otimes l(\sigma)}=
\bigoplus_{\gamma\in I(\lambda)}H_{\gamma}(\lambda),
\end{equation}
where $I(\lambda)\subset (0,1]$ 
denotes the set of eigenvalues of $L_{\sigma\lambda;\lambda}$, and
$H_{\gamma}(\lambda)$ denotes the eigenspace of $L_{\sigma\lambda;\lambda}$ 
corresponding
to the eigenvalue $\gamma\in I(\lambda)$ (we suppress the dependance on $\sigma$
if there is no confusion possible). Observe that $1\in I(\lambda)$ 
and that $L_{\sigma\lambda;\lambda}$ is injective. 

Recall the definition of the set $W^S$ of minimal coset representatives
(cf.\ \eqref{mincoset}). An alternative characterization of $W^S$ is given by
\begin{equation}\label{mincoset2}
W^S=\lbrace \sigma\in W \, | \, \sigma(R_S^+)\subset R^+\rbrace,
\end{equation}
where $R_S^+:=R^+\cap \hbox{span}\lbrace S\rbrace$ 
(cf.\ \cite[Prop.\ 5.1 (iii)]{BGG}).
Using this alternative description of $W^S$ we obtain the following 
properties of $L_{\sigma\lambda;\lambda}$ for $\lambda\in P_{++}(S^c)$.
\begin{Prop}\label{1dimensional}
Suppose that $\sigma\in W^S$ and $\lambda\in P_{++}(S^c)$. 
Then\\ 
{\bf (i)} $L_{\sigma\lambda;\lambda}$ is a compact operator;\\
{\bf (ii)} The eigenspace $H_1(\lambda)$ of $L_{\sigma\lambda;\lambda}$ corresponding
to the eigenvalue $1$ is spanned by the vector $e_0^{\otimes l(\sigma)}$.
\end{Prop}
\begin{proof}
Fix a $\lambda\in P_{++}(S^c)$, and 
let $\sigma=s_{i_1}s_{i_2}\cdots s_{i_l}$ be a reduced expression 
of a minimal coset representative $\sigma\in W^S$.
It is well known that
\begin{equation}\label{constructionbadroots}
R^+\cap \sigma^{-1}(R^-)=\lbrace \gamma_k\rbrace_{k=1}^l,
\end{equation}
where the $\gamma_k$ are defined by \eqref{factohelp}.
We have $\gamma_k\in R^+\setminus R_S^+$ by \eqref{mincoset2}. It follows that
$(\lambda,\gamma_k\spcheck)>0$ for all $k$, since $\lambda\in P_{++}(S^c)$.
By \eqref{action2} and \eqref{facto} it follows that 
$H_1(\lambda)=\hbox{span}\lbrace e_0^{\otimes l(\sigma)}\rbrace$ and that
$H_{\gamma}(\lambda)$ is finite-dimensional for all $\gamma\in
I(\lambda)$.
Since the spectrum of $L_{\sigma\lambda;\lambda}$ (which is equal to 
$I(\lambda)\cup\lbrace 0 \rbrace$) does not have a limit point except
$0$, we conclude that $L_{\sigma\lambda;\lambda}$ is a compact operator 
(cf. \cite[Theorem 12.30]{R}). 
\end{proof}
Let us recall 
the following well known inequalities for weights of finite-dimensional 
irreducible representations of ${\mathfrak{g}}$
(or, equivalently, $U_q(\mathfrak{g})$).
\begin{Prop}\label{needed}
Let $\lambda\in P_+$ and $\mu,\nu\in P(\lambda)$. Then 
$(\lambda,\lambda)\geq (\mu,\nu)$,
and equality holds if and only if $\mu=\nu\in W\lambda$.
\end{Prop}
For a proof of the proposition, see for instance \cite[Prop.\  11.4]{Kac}.
The proof is based on the following lemma, which we will also need
later on. The lemma is a slightly weaker version of \cite[Lemma 11.2]{Kac}.
\begin{Lem}\label{neededlemma}
Let $\lambda\in P_+$ and $\mu\in P(\lambda)\setminus\lbrace \lambda \rbrace$,
and let $m_i\in {\mathbb{Z}}_+$ \textup{(}$i\in [1,r]$\textup{)}  be the
expansion coefficients defined by $\lambda-\mu=\sum_im_i\alpha_i$. 
Then there is an $1\leq i\leq r$ with $m_i>0$ and $\lambda(H_i)\neq 0$.
\end{Lem}
We now have the following proposition, which can be regarded as a
quantum analogue of the ``if'' part of Proposition \ref{semiclassical}.
\begin{Prop}\label{irrep}
Let $\sigma\in W^S$. Then $\pi_\sigma$ restricts to 
an irreducible $*$-represen\-ta\-tion of the factorized $*$-algebra $A_S$. 
In particular, $\pi_\sigma$ restricts to an irreducible $*$-representation
of ${\mathbb{C}}_q\lbrack U/K_S \rbrack$.
\end{Prop}
\begin{proof}
Let $\lambda\in P_{++}(S^c)$ and $\sigma\in W^S$. Suppose 
$H\subset l_2({\mathbb{Z}}_+)^{\otimes l(\sigma)}$ is a non-zero
closed subspace invariant under ${\pi_\sigma}_{\vert A_S}$.
Set $\gamma:=\| {L_{\sigma\lambda;\lambda}}_{\vert H}\|$.
Then $\gamma>0$, since $L_{\sigma\lambda;\lambda}$ is injective
and $\gamma$ is an eigenvalue of $L_{\sigma\lambda;\lambda}/_H$ by
Proposition \ref{1dimensional}{\bf (i)}. 
Let $H_{\gamma}$ be the corresponding eigenspace.
We claim that
\begin{equation}\label{claimbel}
\pi_\sigma((C_{\mu,i;\lambda}^{\lambda})^*C_{\mu,i;\lambda}^{\lambda})
H_{\gamma}=0, \quad  \mu\not=\sigma\lambda.
\end{equation}
Suppose for the moment that the claim is correct. Then 
\eqref{unitarityproperty} and \eqref{claimbel} imply
$\gamma=1$, hence $H_{\gamma}=
\hbox{span}\lbrace e_0^{\otimes l(\sigma)}\rbrace$ by Proposition
\ref{1dimensional}{\bf (ii)}.
So every non-zero closed invariant subspace 
contains the vector $e_0^{\otimes l(\sigma)}$.
Since $H^{\bot}$ is also a closed invariant subspace, we 
must have $H^{\bot}=\lbrace 0 \rbrace$, 
i.e. $H=l_2({\mathbb{Z}}_+)^{\otimes l(\sigma)}$.
Remains therefore to prove the claim \eqref{claimbel}.
By \eqref{part1} we have $\pi_\sigma(C_{\mu,i;\lambda}^{\lambda})=0$
if $\mu<\sigma\lambda$. Hence 
\begin{equation}
\begin{split}
L_{\sigma\lambda;\lambda}
\pi_\sigma((C_{\mu,i;\lambda}^{\lambda})^*C_{\mu,i;\lambda}^{\lambda})&=
q^{(\lambda,\lambda)-(\mu,\sigma\lambda)}
\pi_\sigma((C_{\mu,i;\lambda}^{\lambda}C_{\sigma\lambda;\lambda}^{\lambda})^*
C_{\sigma\lambda;\lambda}^{\lambda}C_{\mu,i;\lambda}^{\lambda})\nonumber\\
&=q^{2(\lambda,\lambda)-2(\mu,\sigma\lambda)}
\pi_\sigma((C_{\mu,i;\lambda}^{\lambda})^*(C_{\sigma\lambda;\lambda}^{\lambda})^*
C_{\sigma\lambda;\lambda}^{\lambda}C_{\mu,i;\lambda}^{\lambda})\nonumber\\
&=q^{2(\lambda,\lambda)-2(\mu,\sigma\lambda)}
\pi_\sigma((C_{\mu,i;\lambda}^{\lambda})^*C_{\sigma\lambda;\lambda}^{\lambda})
\pi_\sigma((C_{\sigma\lambda;\lambda}^{\lambda})^*C_{\mu,i;\lambda}^{\lambda}),
\nonumber
\end{split}
\end{equation}
where we used Proposition \ref{commutation}{\bf (i)} in the second equality
and Proposition \ref{commutation}{\bf (ii)} in the first and third equality.
So \eqref{claimbel} will then follow from
\begin{equation}\label{bijnadaar}
\pi_\sigma((C_{\sigma\lambda;\lambda}^{\lambda})^*
C_{\mu,i;\lambda}^{\lambda})H_{\gamma}= 0, \quad \mu\not=\sigma\lambda,
\end{equation}
in view of the injectivity of $L_{\sigma\lambda;\lambda}$.
Fix $h\in H_{\gamma}$ and $\mu\in P(\lambda)$ with $\mu\not=\sigma\lambda$.
By Lemma \ref{definitionAS} we have 
$(C_{\sigma\lambda;\lambda}^{\lambda})^*C_{\mu;i;\lambda}^{\lambda}\in A_S\subset
{\mathbb{C}}_q\lbrack U/K_S \rbrack$, hence the vector
\begin{equation}\label{image}
\tilde{h}:=\pi_\sigma((C_{\sigma\lambda;\lambda}^{\lambda})^*
C_{\mu;i;\lambda}^{\lambda})h
\end{equation}
lies in the invariant subspace $H$.
Again using the commutation relations given in Proposition \ref{commutation}
and Corollary \ref{other}, we see that $\tilde{h}$  
is an eigenvector of $L_{\sigma\lambda;\lambda}$ with eigenvalue
$\tilde{\gamma}:=q^{2(\lambda,\sigma^{-1}(\mu)-\lambda)}\gamma$. 
We have $\tilde{\gamma}>\gamma$ by Proposition \ref{needed}.
By the maximality of $\gamma$, we conclude that $\tilde{h}=0$.
This proves \eqref{bijnadaar}, hence also the claim \eqref{claimbel}.
\end{proof}
\begin{Def}
We say that the irreducible $*$-representation $\pi_{\sigma}$ 
\textup{(}$\sigma\in W^S$\textup{)} of ${\mathbb{C}}_q\lbrack U/K_S \rbrack$ 
is associated with the Schubert cell $X_{\overline{\sigma}}\subset U/K_S$.
\end{Def} 
The following proposition can be regarded as a quantum analogue
of Proposition \ref{finestructure} as well as of 
the ``only if'' part of Proposition \ref{semiclassical}.
\begin{Prop}\label{notirrep}
Let $\sigma\in W$, and let $\sigma=uv$ be the unique decomposition of $\sigma$ 
with $u\in W^S$ and $v\in W_S$. 
For $\pi_\sigma=\pi_u\otimes\pi_v$ \textup{(}cf.\ \eqref{minimal}\textup{)}
and $t\in T$, we have
\[(\pi_\sigma\otimes\tau_t)(a)=
\pi_u(a)\otimes\hbox{id}^{\otimes l(v)},
\quad a\in {\mathbb{C}}_q\lbrack U/K_S \rbrack.
\]
\end{Prop}
\begin{proof}
Recall that the one-dimensional $*$-representation $\tau_t$ factorizes 
through $\iota_T^*: {\mathbb{C}}_q\lbrack U \rbrack\rightarrow
{\mathbb{C}}_q\lbrack T\rbrack$ 
and that $\pi_i$ factorizes through 
$\phi_i^*: {\mathbb{C}}_q\lbrack U \rbrack\rightarrow 
\mathbb{C}_{q_i}\lbrack SU(2)\rbrack $.
The maps $\iota_T^*$ and $\phi_i^*$ ($i\in S$) factorize through
$\iota_S^*: {\mathbb{C}}_q\lbrack U \rbrack\rightarrow {\mathbb{C}}_q\lbrack
K_S\rbrack$ since the ranges of 
$\iota_T$ and $\phi_i$ ($i\in S$) lie in the Hopf-subalgebra $U_q({\mathfrak{l}}_S)$.
Hence $\pi_v\otimes \tau_t$ ($v\in W_S$, $t\in T$) factorizes through
$\iota_S^*$, say $\pi_v\otimes \tau_t=\pi_{v,t}\circ \iota_S^*$.
Then we have for $a\in {\mathbb{C}}_q\lbrack U/K_S \rbrack$,
\begin{equation}
\begin{split}
(\pi_\sigma\otimes \tau_t)(a)&=
(\pi_u\otimes \pi_v\otimes \tau_t)\circ\Delta(a)\nonumber\\
&=(\pi_u\otimes\pi_{v,t})\circ (\hbox{id}\otimes{\iota}_S^*)\Delta(a)\nonumber\\
&=\pi_u(a)\otimes\pi_{v,t}(1)=\pi_u(a)\otimes\hbox{id}^{\otimes l(v)},\nonumber
\end{split}
\end{equation}
which completes the proof of the proposition.
\end{proof}
\begin{Lem}\label{inequivalent}
The $*$-representations $\lbrace \pi_\sigma\rbrace_{\sigma\in W^S}$,
considered as $*$-representations of $A_S$ 
respectively ${\mathbb{C}}_q\lbrack U/K_S \rbrack$, 
are mutually inequivalent.
\end{Lem}
\begin{proof}
Let $\sigma,\sigma'\in W^S$ with $\sigma\not=\sigma'$ and $\lambda\in P_{++}(S^c)$. 
Then $\sigma\lambda\not=\sigma'\lambda$, since the isotropy subgroup 
$\lbrace \sigma\in W \, | \, 
\sigma\lambda=\lambda\rbrace$ is equal to $W_S$ by Chevalley's Lemma (cf. \cite[Prop.
2.72]{Knapp}).
Without loss of generality we may assume that
$\sigma\lambda\not\geq \sigma'\lambda$. Then we have
$\pi_{\sigma'}((C_{\sigma\lambda,\lambda}^{\lambda})^*C_{\sigma\lambda,\lambda}^{
\lambda})=0$
by \eqref{part1}. On the other hand, $L_{\sigma\lambda;\lambda}$ is injective.
It follows that $\pi_\sigma\not\simeq\pi_{\sigma'}$ as $*$-representations of $A_S$. 
\end{proof}
Let now $\|.\|_u$ be the universal $C^*$-norm 
on ${\mathbb{C}}_q\lbrack U \rbrack$ (cf.\ \cite[\S4]{DK}), so
\begin{equation}
\|a\|_u:=\sup_{\sigma\in W,t\in T}\|(\pi_\sigma\otimes\tau_t)(a)\|,
\quad a\in {\mathbb{C}}_q\lbrack U \rbrack.
\end{equation}
Let $C_q(U)$ (resp. $C_q(U/K_S)$) be the completion of 
${\mathbb{C}}_q\lbrack U \rbrack$ 
(resp.\ ${\mathbb{C}}_q\lbrack U/K_S \rbrack$) 
with respect to $\|.\|_u$. All $*$-representations
$\pi_\sigma\otimes \tau_t$ of ${\mathbb{C}}_q\lbrack U \rbrack$ 
extend to $*$-representations of the $C^*$-algebra $C_q(U)$
by continuity. The results of this section can now be summarized
as follows.
\begin{Thm}\label{restrictionrepr}
Let $S\subsetneq \Delta$. Then 
$\lbrace \pi_\sigma\rbrace_{\sigma\in W^S}$ is a complete set 
of mutually inequivalent irreducible $*$-representations of $C_q(U/K_S)$.
\end{Thm}
\begin{proof}
This follows from the previous results, since every irreducible $*$-represen\-ta\-tion
of $C_q(U/K_S)$  appears as an irreducible component of 
$\sigma_{\vert C_q(U/K_S)}$ 
for some irreducible $*$-represen\-ta\-tion $\sigma$ of $C_q(U)$ 
(cf.\  \cite[Prop.\  2.10.2]{D}).
\end{proof}
Theorem \ref{restrictionrepr} does not imply that 
$\lbrace \pi_\sigma\rbrace_{\sigma\in W^S}$
is a complete set of irreducible $*$-representations of the
$*$-algebra ${\mathbb{C}}_q\lbrack U/K_S \rbrack$ itself.
Indeed, it is not clear that any irreducible $*$-representation
of ${\mathbb{C}}_q\lbrack U/K_S \rbrack$ can be continuously extended to
a $*$-representation of $C_q(U/K_S)$.
In the remainder of this paper we will deal with the classification of the irreducible 
$*$-representations of $A_S$. In particular, this will yield a complete classification
of the irreducible $*$-representations of ${\mathbb{C}}_q\lbrack U/K_S\rbrack$ for the
generalized flag manifolds $U/K_S$ for which the PBW factorization is valid
(cf. Theorem \ref{yes}).
\section{Irreducible $*$-representations of $A_S$}
Let $S\subsetneq \Delta$ be any subset. In this section we show that
$\lbrace \pi_{\sigma}\rbrace_{\sigma\in W^S}$
exhausts the set of irreducible $*$-representations of $A_S$ (up to
equivalence). We fix therefore an arbitrary irreducible $*$-representation
$\tau: A_S\rightarrow {\mathcal{B}}(H)$ and we will show that $\tau\simeq
\pi_{\sigma}$ for a (unique) $\sigma\in W^S$. In order to 
associate the proper minimal coset representative 
$\sigma\in W^S$ with $\tau$, we need to study the range $\tau(A_S)\subset
{\mathcal{B}}(H)$ of $\tau$ in more detail.  
For $\lambda\in P_+(S^c)$ and $\mu,\nu\in P(\lambda)$,
let $\tau^{\lambda}(\mu;\nu),\tau^{\lambda}(\nu)\subset {\mathcal{B}}(H)$ 
be the linear subspaces
\begin{equation}
\begin{split}
\tau^{\lambda}(\mu;\nu) &:=\lbrace \tau((C_{v;v_{\lambda}}^{\lambda})^*
C_{w;v_{\lambda}}^{\lambda}) \, | \, v\in V(\lambda)_{\mu},\,\, w\in V(\lambda)_{\nu}
\rbrace,\\
\tau^{\lambda}(\nu)&:=\lbrace \tau((C_{v;v_{\lambda}}^{\lambda})^*
C_{w;v_{\lambda}}^{\lambda}) \, | \, v\in V(\lambda),\,\, w\in V(\lambda)_{\nu}
\rbrace.
\end{split}
\end{equation}
For $\lambda\in P_+(S^c)$ set 
\begin{equation}\label{D}
D(\lambda):=\lbrace \nu\in P(\lambda) \, | \, \tau^{\lambda}(\nu)\not=\lbrace 0 \rbrace
\rbrace
\end{equation}
and let $D_m(\lambda)$ be the set of weights $\nu\in D(\lambda)$ such that
$\nu'\notin D(\lambda)$ for all $\nu'<\nu$. 
By \eqref{unitarityproperty}, we have $D(\lambda)\not=\emptyset$, hence
also $D_m(\lambda)\not=\emptyset$.
We start with a lemma which is useful for the computation of commutation relations
in $\tau(A_S)\subset {\mathcal{B}}(H)$.
\begin{Lem}\label{help!}
Let $\lambda,\Lambda\in P_+(S^c)$ and $\nu\in D_m(\lambda)$. 
Let $v\in V(\lambda)$, $v'\in V(\lambda)_{\nu'}$ with $\nu'<\nu$ 
and $w,w'\in V(\Lambda)$.
Then the product of the four matrix elements 
$(C_{v;v_{\lambda}}^{\lambda})^*$, $C_{v';v_{\lambda}}^{\lambda}$,
$(C_{w;v_{\Lambda}}^{\Lambda})^*$, and  $C_{w';v_{\Lambda}}^{\Lambda}$,
taken in an arbitary order, is contained in $\hbox{Ker}(\tau)$.
\end{Lem}
\begin{proof}
Since $\hbox{Ker}(\tau)$ is a two-sided $*$-ideal in $A_S$, it follows
from the definitions that
\[(C_{w;v_{\Lambda}}^{\Lambda})^*C_{w';v_{\Lambda}}^{\Lambda}
(C_{v;v_{\lambda}}^{\lambda})^*C_{v';v_{\lambda}}^{\lambda}\in \hbox{Ker}(\tau).\]
If the product of the four matrix coefficients is taken in a different order,
then we can rewrite it by Proposition \ref{commutation} and 
by Corollary \ref{other} as a linear combination of products of matrix elements
\[(C_{u;v_{\Lambda}}^{\Lambda})^*C_{u';v_{\Lambda}}^{\Lambda}
(C_{x;v_{\lambda}}^{\lambda})^*C_{x';v_{\lambda}}^{\lambda}\]
with $x'\in V(\lambda)_{\nu''}$ and $\nu''\leq\nu'<\nu$. These 
are all contained in $\hbox{Ker}(\tau)$, since $\nu\in D_m(\lambda)$.
\end{proof}
\begin{Lem}\label{w}
Let $\lambda\in P_+(S^c)$ and $\nu\in D_m(\lambda)$. Then\\ 
{\bf (i)} $\tau^{\lambda}(\nu;\nu)\not=\lbrace 0 \rbrace$;\\
{\bf (ii)} $\nu=\sigma\lambda$ for some $\sigma\in W^S$.
\end{Lem}
\begin{proof}
Let $\lambda\in P_+(S^c)$ and $\nu\in D_m(\lambda)$. 
Fix weight vectors $v\in V(\lambda)_{\mu}$, $w\in V(\lambda)_{\nu}$ such that 
$T_{v;w}:=\tau((C_{v;v_{\lambda}}^{\lambda})^*C_{w;v_{\lambda}}^{\lambda})\not=0$.
By Lemma \ref{help!}, we compute
\begin{equation}\label{subres}
\begin{split}
(T_{v;w})^*T_{v;w}&=q^{(\mu,\nu)-(\lambda,\lambda)}
\tau(C_{v;v_{\lambda}}^{\lambda}(C_{v;v_{\lambda}}^{\lambda}
C_{w;v_{\lambda}}^{\lambda})^*C_{w;v_{\lambda}}^{\lambda})\\
&=\tau(C_{v;v_{\lambda}}^{\lambda}(C_{v;v_{\lambda}}^{\lambda})^*)
T_{w;w},
\end{split}
\end{equation}
where we used Proposition \ref{commutation}{\bf (ii)} in the first equality
and Proposition \ref{commutation}{\bf (i)} in the second equality.
On the other hand, $(T_{v;w})^*T_{v;w}\not=0$ since ${\mathcal{B}}(H)$ is a
$C^*$-algebra, so we conclude that $T_{w;w}\not=0$. In particular, 
$\tau^{\lambda}(\nu,\nu)\not=\lbrace 0 \rbrace$. Formula \eqref{subres}
for $v=w$ gives
\[
0\not=(T_{w;w})^*T_{w;w}=
\tau(C_{w;v_{\lambda}}^{\lambda}(C_{w;v_{\lambda}}^{\lambda})^*)T_{w;w}
=q^{(\lambda,\lambda)-(\nu,\nu)}T_{w;w}T_{w;w},
\]
where we have used Proposition \ref{commutation}{\bf (ii)} in the last
equality. It follows that $(\lambda,\lambda)=(\nu,\nu)$, since $T_{w;w}$ is
selfadjoint.
By Proposition \ref{needed} 
we obtain $\nu=\sigma\lambda$ for some $\sigma\in W^S$. 
\end{proof}
For $\lambda\in P_+(S^c)$ and $\nu\in D_m(\lambda)$ we set
\begin{equation}\label{Lw}
L_{\nu;\lambda}:=\tau((C_{\nu;\lambda}^{\lambda})^*
C_{\nu;\lambda}^{\lambda}).
\end{equation}
This definition makes sense since $\hbox{dim}(V(\lambda)_{\nu})=1$ by Lemma
\ref{w}{\bf (ii)}.
Furthermore, $L_{\nu;\lambda}$ is a non-zero selfadjoint operator which commutes 
with the elements of $\tau(A_S)$ in the following way.
\begin{Lem}\label{derivedcomm}
Let $\lambda,\Lambda\in P_+(S^c)$ and $\nu\in D_m(\lambda)$.
For $v\in V(\Lambda)_{\mu}$, $w\in V(\Lambda)_{\mu'}$ we have
\[
L_{\nu;\lambda}\tau((C_{v;v_{\Lambda}}^{\Lambda})^*C_{w;v_{\Lambda}}^{\Lambda})=
q^{2(\nu,\mu'-\mu)}
\tau((C_{v;v_{\Lambda}}^{\Lambda})^*C_{w;v_{\Lambda}}^{\Lambda})L_{\nu;\lambda}.
\]
\end{Lem}
\begin{proof}
By Lemma \ref{help!} and the commutation relations in section 3 we compute
\begin{equation}
\begin{split}
L_{\nu;\lambda}\tau((C_{v;v_{\Lambda}}^{\Lambda})^*C_{w;v_{\Lambda}}^{\Lambda})&=
q^{(\lambda,\Lambda)-(\nu,\mu)}\tau((C_{v;v_{\Lambda}}^{\Lambda}
C_{v_{\nu};v_{\lambda}}^{\lambda})^*C_{v_{\nu};v_{\lambda}}^{\lambda}
C_{w;v_{\Lambda}}^{\Lambda})\nonumber\\
&=q^{2(\lambda,\Lambda)-2(\nu,\mu)}\tau((C_{v;v_{\Lambda}}^{\Lambda})^*
(C_{v_{\nu};v_{\lambda}}^{\lambda})^*C_{v_{\nu};v_{\lambda}}^{\lambda}
C_{w;v_{\Lambda}}^{\Lambda})\nonumber\\
&=q^{(\nu,\mu')+(\lambda,\Lambda)-2(\nu,\mu)}\tau((C_{v;v_{\Lambda}}^{\Lambda})^*
(C_{v_{\nu};v_{\lambda}}^{\lambda})^*C_{w;v_{\Lambda}}^{\Lambda}
C_{v_{\nu};v_{\lambda}}^{\lambda})\nonumber\\
&=q^{2(\nu,\mu'-\mu)}
\tau((C_{v;v_{\Lambda}}^{\Lambda})^*C_{w;v_{\Lambda}}^{\Lambda})L_{\nu;\lambda},
\nonumber
\end{split}
\end{equation}
where we used Proposition \ref{commutation}{\bf (ii)} for the first and
fourth equality, Proposition \ref{commutation}{\bf (i)} for the second
equality, and Corollary \ref{other} for the third equality. 
\end{proof}
It follows from Lemma \ref{derivedcomm} that $\hbox{Ker}(L_{\nu;\lambda})\subsetneq H$ 
is a closed invariant subspace. By the irreducibility of $\tau$, we thus obtain the
following corollary.
\begin{Cor}\label{injective}
Let $\lambda\in P_+(S^c)$ and $\nu\in D_m(\lambda)$. Then
$L_{\nu;\lambda}$ is injective.
\end{Cor}
The minimal coset representative $\sigma$ of Lemma \ref{w}{\bf (ii)}
is unique and independent of $\lambda\in P_+(S^c)$ in the following sense.
\begin{Lem}\label{uniqueness}
There exists a unique $\sigma\in W^S$ 
such that $D_m(\lambda)=\lbrace \sigma\lambda\rbrace$ 
for all $\lambda\in P_+(S^c)$. 
\end{Lem}
\begin{proof}
Let $\Lambda\in P_{++}(S^c)$ and $\nu\in D_m(\Lambda)$. Then there exists
a unique $\sigma\in W^S$ such that $\nu=\sigma\Lambda$ by Lemma \ref{w}{\bf (ii)}
and by Chevalley's Lemma (cf. \cite[Prop.\ 2.27]{Knapp}).
Fix furthermore arbitrary $\lambda\in P_+(S^c)$ and 
$\nu'\in D_m(\lambda)$. Choose a $\sigma'\in W$ such that $\nu'=\sigma'\lambda$.
By Lemma \ref{help!} and the commutation relations of section 3, we compute
\begin{equation}
\begin{split}
L_{\nu;\Lambda}L_{\nu';\lambda}&=q^{(\Lambda,\lambda)-(\nu,\nu')}
\tau((C_{\nu';\lambda}^{\lambda}C_{\nu;\Lambda}^{\Lambda})^*
C_{\nu;\Lambda}^{\Lambda}C_{\nu';\lambda}^{\lambda})\nonumber\\
&=q^{3(\Lambda,\lambda)-3(\nu,\nu')}\tau(
(C_{\nu';\lambda}^{\lambda})^*(C_{\nu;\Lambda}^{\Lambda})^*
C_{\nu';\lambda}^{\lambda}C_{\nu;\Lambda}^{\Lambda})\nonumber\\
&=q^{2(\Lambda,\lambda)-2(\nu,\nu')}L_{\nu';\lambda}L_{\nu;\Lambda},\nonumber
\end{split}
\end{equation}
where we used Proposition \ref{commutation}{\bf (ii)} in the first
and third equality and Proposition \ref{commutation}{\bf (i)} twice in the
second equality.
If we repeat the same computation, but now using Corollary
\ref{other} twice in the second equality, then we obtain
\[
L_{\nu;\Lambda}L_{\nu';\lambda}=
q^{2(\nu,\nu')-2(\Lambda,\lambda)}L_{\nu';\lambda}L_{\nu;\Lambda},
\]
hence
\[\bigl(q^{2(\Lambda,\lambda)-2(\nu,\nu')}-q^{2(\nu,\nu')-2(\Lambda,\lambda)}
\bigr)L_{\nu';\lambda}L_{\nu;\Lambda}=0.\]
By Corollary \ref{injective} we have $L_{\nu';\lambda}L_{\nu;\Lambda}\not=0$, so we
conclude that
\[
(\Lambda,\lambda)-(\nu,\nu')=(\Lambda,\lambda-\sigma^{-1}\sigma'\lambda)=0.
\] 
Since $\Lambda\in P_{++}(S^c)$ and $\lambda\in P_+(S^c)$, 
it follows from Lemma \ref{neededlemma}
that $\lambda=\sigma^{-1}\sigma'\lambda$, i.e. 
$\nu'=\sigma\lambda$. Hence, $D_m(\lambda)=\lbrace
\sigma\lambda\rbrace$ for all $\lambda\in P_+(S^c)$. 
\end{proof}
In the remainder of this section we write $\sigma$ for 
the unique minimal coset representative such that $D_m(\lambda)=\lbrace
\sigma\lambda\rbrace$ for all $\lambda\in P_+(S^c)$. 
We are going to prove that $\tau\simeq\pi_{\sigma}$.
First we look for the analogue of the distinguished vector
$e_0^{\otimes l(\sigma)}$ (cf.\  Proposition \ref{1dimensional}{\bf (ii)})
in the representation space $H$ of $\tau$.

The spectrum $I(\lambda)$ of $L_{\sigma\lambda;\lambda}$ 
is contained in $[0,\infty)$, since
$L_{\sigma\lambda;\lambda}$ is a positive operator. 
By considering the spectral decomposition
of $L_{\sigma\lambda;\lambda}$, one obtains the following
corollary of Lemma \ref{derivedcomm} and \cite[Lemma 4.3]{Ko}. 
\begin{Cor}\label{spec}
Let $\lambda\in P_+(S^c)$. Then
$I(\lambda)\subset [0,\infty)$
is a countable set with no limit points, except possibly $0$.
\end{Cor}
The proof of Corollary \ref{spec} is similar to the proof of
\cite[Prop.\ 3.9]{S} and of \cite[Prop.\ 4.2]{Ko}.

By Corollary \ref{spec} we have an orthogonal direct sum decomposition
\begin{equation}\label{eig}
H=\bigoplus_{\gamma\in I(\lambda)\cap {\mathbb{R}}_{>0}}H_{\gamma}(\lambda)
\end{equation}
into eigenspaces of $L_{\sigma\lambda;\lambda}$, where $H_{\gamma}(\lambda)$
is the eigenspace of $L_{\sigma\lambda;\lambda}$ 
corresponding to the eigenvalue $\gamma$.
Let $\gamma_0(\lambda)>0$ be the largest eigenvalue of $L_{\sigma\lambda;\lambda}$.
\begin{Lem}\label{zero}
Let $\lambda\in P_+(S^c)$, $v\in V(\lambda)$, $w\in V(\lambda)_{\nu}$
and assume that $\nu\not=\sigma\lambda$. Then
$\tau((C_{v;v_{\lambda}}^{\lambda})^*C_{w;v_{\lambda}}^{\lambda})
(H_{\gamma_0(\lambda)}(\lambda))=
\lbrace 0 \rbrace$.
\end{Lem}
\begin{proof}
Let $\lambda\in P_+(S^c)$, $v\in V(\lambda)_{\mu}$ and $w\in V(\lambda)_{\nu}$. 
By Lemma \ref{help!} and the commutation relations in section 3, we
compute
\begin{equation}
\begin{split}
L_{\sigma\lambda;\lambda}\tau((C_{v;v_{\lambda}}^{\lambda})^*C_{w;v_{\lambda}}^{
\lambda})&=\tau(C_{v_{\sigma\lambda};v_{\lambda}}^{\lambda}
(C_{v;v_{\lambda}}^{\lambda}C_{v_{\sigma\lambda};v_{\lambda}}^{\lambda})^*
C_{w;v_{\lambda}}^{\lambda})\nonumber\\
&=q^{(\lambda,\lambda)-(\mu,\sigma\lambda)}
\tau(C_{v_{\sigma\lambda};v_{\lambda}}^{\lambda}
(C_{v;v_{\lambda}}^{\lambda})^*(C_{v_{\sigma\lambda};v_{\lambda}}^{\lambda})^*
C_{w;v_{\lambda}}^{\lambda})\nonumber\\
&=q^{2(\lambda,\lambda)-2(\mu,\sigma\lambda)}
\tau((C_{v;v_{\lambda}}^{\lambda})^*C_{v_{\sigma\lambda};v_{\lambda}}^{\lambda})
\tau((C_{v_{\sigma\lambda};v_{\lambda}}^{\lambda})^*C_{w;v_{\lambda}}^{\lambda}),
\nonumber
\end{split}
\end{equation}
where we used Proposition \ref{commutation}{\bf (i)} in the second equality
and Proposition \ref{commutation}{\bf (ii)} in the first and third equality.
This computation, together with the injectivity of $L_{\sigma\lambda;\lambda}$,
shows that it suffices to give a proof of the lemma for the special case
that $v=v_{\sigma\lambda}$.
So we fix $h\in H_{\gamma_0(\lambda)}(\lambda)$ and 
$w\in V(\lambda)_{\nu}$ with $\nu\in P(\lambda)$ and $\nu\not=\sigma\lambda$. 
It follows from Lemma \ref{derivedcomm} that
$\tilde{h}:=
\tau((C_{v_{\sigma\lambda};v_{\lambda}}^{\lambda})^*C_{w;v_{\lambda}}^{\lambda})h$
is an eigenvector of $L_{\sigma\lambda;\lambda}$ with eigenvalue 
$\tilde{\gamma}_0(\lambda)=q^{2(\lambda,\sigma^{-1}(\nu)-\lambda)}\gamma_0(\lambda)$.
By Proposition \ref{needed} we have $\tilde{\gamma}_0(\lambda)>\gamma_0(\lambda)$, 
hence $\tilde{h}=0$ by the maximality of the eigenvalue $\gamma_0(\lambda)$.
\end{proof}
\begin{Cor}\label{Cor2}
$\gamma_0(\lambda)=1$ for all $\lambda\in P_+(S^c)$.
\end{Cor}
\begin{proof}
Follows from \eqref{unitarityproperty} and Lemma \ref{zero}.
\end{proof} 
The linear subspace of ${\mathbb{C}}_q\lbrack U \rbrack$ spanned by the matrix elements 
$\lbrace C_{\sigma\mu;\mu}^{\mu}\rbrace_{\mu\in P_+}$
is a subalgebra of ${\mathbb{C}}_q\lbrack U \rbrack$ with algebraic generators
$C_{\sigma\varpi_i;\varpi_i}^{\varpi_i}$ ($i\in [1,r]$), since 
$C_{\sigma\mu;\mu}^{\mu}C_{\sigma\nu;\nu}^{\nu}=
\lambda_{\mu,\nu}C_{\sigma(\mu+\nu);\mu+\nu}^{\mu+\nu}$,
where the scalar $\lambda_{\mu,\nu}\in {\mathbb{T}}$
depends on the particular choices of orthonormal bases for the finite-dimensional 
irreducible representations $V(\mu)$ and $V(\nu)$
(cf. \cite[Proof of Prop.\ 3.12]{S}).
Then it follows from Proposition \ref{commutation} and Lemma \ref{help!} that 
\begin{equation}\label{product}
L_{\sigma(\mu+\nu);\mu+\nu}=L_{\sigma\mu;\mu}L_{\sigma\nu;\nu}
\end{equation}
for all $\mu,\nu\in P_+(S^c)$,
hence $\hbox{span}\lbrace L_{\sigma\lambda;\lambda}\rbrace_{\lambda\in P_+(S^c)}$ 
is a commutative subalgebra of ${\mathcal{B}}(H)$. 
Set
\begin{equation}
H_1:=\bigcap_{i\in S^c}H_1(\varpi_i),
\end{equation}
then $H_1\subset H_1(\lambda)$ for all $\lambda\in P_+(S^c)$ by \eqref{product}.
\begin{Lem}
$H_1=H_1(\lambda)$ for all $\lambda\in P_{++}(S^c)$.
In particular, $H_1\not=\lbrace 0 \rbrace$.
\end{Lem}
\begin{proof}
For $\mu\in P_+(S^c)$ we have $\|L_{\sigma\mu;\mu}\|=1$. Moreover,
for any $h\in H$, 
\begin{equation}\label{equiva}
h\in H_1(\mu) \quad \Leftrightarrow \quad \|L_{\sigma\mu;\mu}h\|=\|h\|.
\end{equation}
This follows from the
eigenspace decomposition \eqref{eig} for $L_{\sigma\mu;\mu}$ 
and the fact that $1$ is the
largest eigenvalue of $L_{\sigma\mu;\mu}$.
Let $\lambda\in P_{++}(S^c)$ and choose arbitrary $i\in S^c$. Then
$\lambda=\mu+\varpi_i$ for certain $\mu\in P_+(S^c)$. By \eqref{product}, 
we obtain for $h\in H_1(\lambda)$,
\begin{equation}
\begin{split}
\|h\|=\|L_{\sigma\lambda;\lambda}h\|&=
\|L_{\sigma\mu;\mu}L_{\sigma\varpi_i;\varpi_i}h\|\nonumber\\
&\leq \|L_{\sigma\varpi_i;\varpi_i}h\|\leq\|h\|,\nonumber
\end{split}
\end{equation}
hence we have equality everywhere. By \eqref{equiva}, it follows that $h\in
H_1(\varpi_i)$.
Since $i\in S^c$ was arbitrary, we conclude that $h\in H_1$.
\end{proof}
\begin{Lem}\label{tussenstukje}
Let $\lambda\in P_+(S^c)$. For all $v\in V(\lambda)_{\mu}$ with
$\mu\not=\sigma\lambda$ we have
\[\tau\bigl((C_{v;v_{\lambda}}^{\lambda})^*
C_{v_{\sigma\lambda};v_{\lambda}}^{\lambda}\bigr)(H_1)\subset H_1^{\bot}.\]
\end{Lem}
\begin{proof}
Let $\Lambda\in P_{++}(S^c)$, $\lambda\in P_+(S^c)$, and $v\in V(\lambda)_{\mu}$ with
$\mu\not=\sigma\lambda$ and $\mu\in P(\lambda)$. Then 
\begin{equation}
L_{\sigma\Lambda;\Lambda}\tau((C_{v;v_{\lambda}}^{\lambda})^*
C_{v_{\sigma\lambda};v_{\lambda}}^{\lambda})=q^{2(\Lambda,\lambda-\sigma^{-1}(\mu))}
\tau((C_{v;v_{\lambda}}^{\lambda})^*
C_{v_{\sigma\lambda};v_{\lambda}}^{\lambda})L_{\sigma\Lambda;\Lambda}
\end{equation}
by Lemma \ref{derivedcomm}. By Lemma \ref{neededlemma} we have
$(\Lambda,\lambda-\sigma^{-1}(\mu))>0$. Hence, 
\begin{equation}
\begin{split}
\tau((C_{v;v_{\lambda}}^{\lambda})^*
C_{v_{\sigma\lambda};v_{\lambda}}^{\lambda})(H_1)
&=\tau((C_{v;v_{\lambda}}^{\lambda})^*
C_{v_{\sigma\lambda};v_{\lambda}}^{\lambda})(H_1(\Lambda))\nonumber\\
&\subset \bigoplus_{\gamma<1}H_{\gamma}(\Lambda)
=H_1(\Lambda)^{\bot}=H_1^{\bot},\nonumber
\end{split}
\end{equation}
which completes the proof of the lemma.
\end{proof}
\begin{Cor}\label{1}
$\dim(H_1)=1$.
\end{Cor}
\begin{proof}
By Lemma \ref{zero} and Lemma \ref{tussenstukje} we obtain for any $0\not=h\in H_1$,
\[
\overline{\tau(A_S)h}\subset \hbox{span}\lbrace h \rbrace \oplus H_1^{\bot},
\]
where the overbar means closure.
By the irreducibility of $\tau$, we conclude that
$\hbox{span}\lbrace h \rbrace=H_1$. 
\end{proof}
Any vector $h\in H_1$ with $\|h\|=1$ can serve now as 
the analogue in the representation space $H$ of 
the distinguished vector $e_0^{\otimes l(\sigma)}$ in the representation
space of $\pi_{\sigma}$. By comparing the Gel'fand-Naimark-Segal states
of $\tau$ and $\pi_{\sigma}$ taken with respect to the cyclic vector
$h\in H_1$ ($\|h\|=1$) resp.\  $e_0^{\otimes l(\sigma)}$, we obtain
the following lemma.
\begin{Lem}
We have $\tau\simeq\pi_\sigma$ as irreducible $*$-representations 
of $A_S$.
\end{Lem}
\begin{proof}
Fix an $h\in H_1$ with $\|h\|=1$, and define the Gel'fand-Naimark-Segal states
$\phi_{\tau},\phi_{\pi_\sigma}:A_S\rightarrow {\mathbb{C}}$ by
\begin{equation}
\phi_{\tau}(a):=(\tau(a)h,h),\,\quad 
\phi_{\pi_\sigma}(a):=(\pi_\sigma(a)e_0^{\otimes l(\sigma)},e_0^{\otimes l(\sigma)}).
\end{equation}
Then we have for $\phi=\phi_{\tau}$ (resp. $\phi=\phi_{\pi_\sigma}$),
\begin{equation}\label{GNS}
\phi((C_{\mu,i;\lambda}^{\lambda})^*C_{\nu,j;\lambda}^{\lambda})
=\delta_{\mu,\sigma\lambda}\delta_{\nu,\sigma\lambda}
\end{equation}
for $\lambda\in P_+(S^c)$, 
$\mu,\nu\in P(\lambda)$, $i\in [1,\hbox{dim}(V(\lambda)_{\mu})]$,
and $j\in [1,\hbox{dim}(V(\lambda)_{\nu})]$. Indeed,
\eqref{GNS} for $\phi=\phi_{\tau}$ follows from Lemma \ref{zero} and Lemma
\ref{tussenstukje}.
For $\phi=\phi_{\pi_\sigma}$, recall that $\pi_\sigma$ is an irreducible 
$*$-representation of $A_S$ (Proposition \ref{irrep}). 
We have seen in the previous section
that $L_{\sigma\lambda;\lambda}=\pi_{\sigma}((C_{\sigma\lambda;\lambda}^{\lambda})^* 
C_{\sigma\lambda;\lambda}^{\lambda})$ is injective for all $\lambda\in P_+(S^c)$,
hence $\sigma\lambda\in D(\lambda)$ (cf.\ \eqref{D}) 
for all $\lambda\in P_+(S^c)$. By \eqref{part1},
we actually have $\sigma\lambda\in D_m(\lambda)$ for all $\lambda\in P_+(S^c)$. 
Hence the labeling $\sigma\in W^S$
of $\pi_{\sigma}$ coincides with its (unique) minimal coset representative 
defined in Lemma \ref{uniqueness}.
Furthermore, the one-dimensional subspace
$H_1$ for $\pi_\sigma$ is equal to $\hbox{span}\lbrace e_0^{\otimes l(\sigma)}
\rbrace$ (cf. Proposition \ref{1dimensional}{\bf (ii)}, Lemma \ref{1}). 
So \eqref{GNS} for $\phi=\phi_{\pi_{\sigma}}$ follows 
again from Lemma \ref{zero} and Lemma \ref{tussenstukje}.

By linearity it follows from \eqref{GNS} that $\phi_{\tau}=\phi_{\pi_\sigma}$, hence
$\tau$ and $\pi_\sigma$ are unitarily equivalent $*$-representations
(cf. \cite[Prop.\ 2.4.1]{D}).
\end{proof}
We may summarize the results of this section as follows.
\begin{Thm}\label{lastjob}
For all $S\subsetneq \Delta$,
$\lbrace \pi_\sigma\rbrace_{\sigma\in W^S}$ is a complete set of mutually inequivalent,
irreducible $*$-representations of the factorized $*$-subalgebra $A_S$.
\end{Thm}
Combining Proposition \ref{notirrep}, Theorem \ref{yes} and Theorem \ref{lastjob}
we obtain the following theorem. 
\begin{Thm}\label{last}
$\lbrace \pi_\sigma\rbrace_{\sigma\in W^S}$ is a complete set of mutually inequivalent,
irreducible $*$-representations of ${\mathbb{C}}_q\lbrack U/K_S \rbrack$ in the
the following cases:\\
{\bf (i)} $S=\emptyset$, i.e. $U/K_S=U/T$ is the full flag manifold;\\
{\bf (ii)} $\#S^c=1$ and the simple root $\alpha\in S^c$ is a Gel'fand node.\\
For these cases the restriction to ${\mathbb{C}}_q[U/K_S]$ of the universal 
$C^*$-norm on ${\mathbb{C}}_q\lbrack U \rbrack$ coincides with the universal
$C^*$-norm on ${\mathbb{C}}_q\lbrack U/K_S \rbrack$.
\end{Thm}
\section*{Acknowledgements}
The research for this work was started while the first author was a guest 
at the University of Kobe during a period of ten weeks in the Spring of 1997,
for which he received financial support by NWO/Nissan.
He greatly acknowledges the hospitality of the Department of Mathematics 
at the University of Kobe. The authors would like to thank Erik T. Koelink 
for helpful discussions about the topic of this paper. The authors thank
Professor P. Littelmann for pointing out the connection of the algebra $A_S^0$ 
with doubled $G$-varieties (cf. Remark \ref{Litremark}).

\bibliographystyle{amsplain}

\end{document}